\documentclass[11pt,a4paper]{article}
\usepackage{soul}
\usepackage{amssymb}
\usepackage{amsmath}
\usepackage{latexsym}
\usepackage{graphicx} 
\usepackage{hyperref}
\usepackage{amsthm}
\usepackage{verbatim}
\usepackage{psfrag}
\usepackage{bbm}
\usepackage{dsfont}
\usepackage{mathrsfs}
\usepackage{tikz,pgfplots}
\usepackage{yfonts}
\usepackage{color}
\usepackage{enumerate}
\usepackage{subcaption}
\usepackage{caption}
\usepackage{pgfplots}
\usepackage{subcaption}
\usepackage{multirow}
\usepackage{algpseudocode}
\usepackage{algorithm}

\newcommand*{\collauthor}[2]{{#1}$^{#2}$}

\newcommand*{\affiliation}[2]{$\mbox{}^{{#2}}${#1}}
\newcommand*{\colltitle}[1]{\textbf{#1}}

\newtheorem{pro}{Proposition}[section]

\newtheorem{remark}[pro]{Remark}

\newtheorem{framework}{Framework}[section]
\newenvironment{keywords}[1]{\vspace{1cm}\\{\bf \slshape{Keywords}}\quad\slshape{#1}}{}

\DeclareMathOperator*{\argmin}{argmin}

\theoremstyle{definition}
\newtheorem{exmp}{Example}[]
\DeclareMathAlphabet\mathbfcal{OMS}{cmsy}{b}{n}

\begin{document}
\begin{center}
\begin{Large}
  \colltitle{{Deep learning algorithms for solving high dimensional nonlinear backward stochastic differential equations}}
\end{Large}
\vspace*{1.5ex}

\begin{sc}
\begin{large}
\collauthor{Lorenc Kapllani}{},
\collauthor{Long Teng}{}
\end{large}
\end{sc}
\vspace{1.5ex}

\affiliation{Lehrstuhl f\"ur Angewandte Mathematik und Numerische Analysis,\\
Fakult\"at f\"ur Mathematik und Naturwissenschaften,\\
Bergische Universit\"at Wuppertal, Gau{\ss}str. 20, \\
42119 Wuppertal, Germany\linebreak }{} \\

\end{center}

\section*{Abstract}
 In this work, we propose a new deep learning-based scheme for solving high dimensional nonlinear backward stochastic differential equations (BSDEs). The idea is to reformulate the problem as a global optimization, where the local loss functions are included. Essentially, we approximate the unknown solution of a BSDE using a deep neural network and its gradient with automatic differentiation. The approximations are performed by globally minimizing the quadratic local loss function defined at each time step, which always includes the terminal condition. This kind of loss functions are obtained by iterating the Euler discretization of the time integrals with the terminal condition. Our formulation can prompt the stochastic gradient descent algorithm not only to take the accuracy at each time layer into account, but also converge to a good local minima. In order to demonstrate performances of our algorithm, several high-dimensional nonlinear BSDEs including pricing problems in finance are provided.
\begin{keywords}
Backward stochastic differential equations, High dimensional problems, deep neural network, recurrent neural network, automatic differentiation, Iterative discretization, Nonlinear option pricing
\end{keywords}

\section{Introduction}
\label{sec1}
In this work we consider the decoupled forward backward stochastic differential equation (BSDE) of the form
\begin{equation}
    \begin{split}
        \left\{
            \begin{array}{rcl}
                dX_t & = & \mu \left(t, X_t\right)\,dt + \sigma \left(t, X_t\right)\,dW_t, \quad X_0 = x_0,\\
   	   		    -dY_t & = & f\left(t, X_t, Y_t, Z_t\right)\,dt -Z_t\,dW_t,\\  
   	   	    	Y_T & = & \xi = g\left(X_T\right),
            \end{array}
        \right. 
    \end{split}
\label{eq1}
\end{equation}
where $X_t, \mu \in \mathbb{R}^d $, $\sigma$ is a $d\times d$ matrix, $W_t = \left( W_t^1, \cdots, W_t^d \right)^\top$ is a $d$-dimensional Brownian motion, $f\left(t,X_t,Y_t,Z_t\right):\left[0,T\right]\times\mathbb{R}^d\times\mathbb{R}\times\mathbb{R}^{1\times d} \to\mathbb{R}$ is the driver function and $\xi$ is the terminal condition which depends on the final value of the forward 
stochastic differential equation (SDE), $X_T.$ The existence and uniqueness of the solution of~\eqref{eq1} are proven in ~\cite{Pardoux1990}. After that, BSDEs have found various applications in finance. For example, as the first claim of applications in finance, it has been shown in ~\cite{El1997} that the price and delta hedging of an option can be represented as a BSDE, and many others such as jump-diffusion models~\cite{eyraud2005backward}, defaultable options~\cite{ankirchner2010credit}, local volatility models~\cite{labart2011parallel}, stochastic volatility models~\cite{fahim2011probabilistic}.

In most cases BSDEs cannot be solved explicitly, advanced numerical techniques to approximate BSDE solutions become desired, especially for the high-dimensional nonlinear BSDEs. In the recent years, many various numerical methods have been proposed for solving BSDEs, e.g., \cite{bouchard2004discrete,zhang2004numerical,gobet2005regression,lemor2006rate,zhao2006new,bender2008,ma2008numerical,zhao2010stable,gobet2010solving,crisan2012solving,zhao2014new,ruijter2015fourier,ruijter2016fourier} and many others. However, most of them are not suitable for solving high-dimensional BSDEs due to the exponentially increasing computational cost with the dimensionality. Although some methods or techniques can be used to accelerate the computations, e.g., methods on sparse grids or parallel computations in graphics processing unit (GPU), only the moderate dimensional BSDEs can be solved numerically for reasonable computational time. We refer to, e.g., \cite{zhang2013sparse,fu2017efficient,chassagneux2021learning} for the methods on sparse grids, and \cite{gobet2016stratified, kapllani2022multistep} for the GPU-based parallel computing.

Recently, several different types of approaches have been proposed to solve high dimensional BSDEs: the multilevel Monte Carlo method based on Picard iteration~\cite{weinan2019multilevel,becker2020numerical,hutzenthaler2021overcoming,hutzenthaler2021multilevel,nguyen2022multilevel}; the regression tree-based methods ~\cite{teng2021review,teng2022gradient}; deep learning-based methods~\cite{weinan2017deep,han2018solving,raissi2018forward,wang2018deep,pereira2019learning,fujii2019asymptotic,ji2020three,hure2020deep,gnoatto2020deep,kremsner2020deep,beck2021deep,chen2021deep,jiang2021convergence,liang2021deep,ji2021control,negyesi2021one,pham2021neural,takahashi2022new,germain2022approximation,andersson2022convergence,ji2022deep}. The first deep learning-based method to approximate the high dimensional BSDEs was proposed in \cite{weinan2017deep}, which has been extended and further studied. It has been pointed out that the method in \cite{weinan2017deep} suffers from the following demerits: 1. It can be stuck in poor local minima or even diverge, especially for a complex solution structure and a long terminal time, see, e.g., \cite{hure2020deep}. 2. It is only capable of achieving good approximations of $Y_0$ and $Z_0,$ namely the solution of a BSDE at the initial time, see \cite{raissi2018forward}.

Motivated by the demerits above we present a novel deep learning-based algorithm to approximate the solutions of nonlinear high dimensional BSDEs. The essential concept is to formulate the problem as a global optimization with local loss functions including the terminal condition. Our formulation is obtained by using the Euler discretization of the time integrals and iterating it with the terminal condition, i.e., iterative time discretization, this might be seen also as a multi-step time discretization. The algorithm estimates the unknown solution (the $Y$ process) using a deep neural network  and its gradient (the $Z$ process) via automatic differentiation (AD). These approximations are performed from the global minimization of the local loss functions defined at each time point from the iterative time discretization. In~\cite{raissi2018forward}, the author have introduced a similar strategy based on local loss functions arising from Euler discretization at each time interval, with the terminal condition included as an additional term in the loss function, i.e., the proposed algorithm attempt to match the dynamics of the BSDE at each time interval. This approach achieves a good approximation of processes $Y$ and $Z$ not only at the initial time but also at each time layer. Hence, it can overcome the second demerit in~\cite{weinan2017deep}. However, the scheme in~\cite{raissi2018forward} still suffers for the first demerit, it can be stuck in poor local minima for the problems with a highly complex structure and a long terminal time, this will be demonstrated in our numerical experiments in Sec. \ref{sec4}. Note that it does not help the stochastic gradient descent (SGD) algorithm in~\cite{weinan2017deep,raissi2018forward} to converge to a good local minima just by considering another network architecture. For instance, the recurrent neural network (RNN) type architectures are specialized for learning long complex sequences. However, it has been pointed out in~\cite{hure2020deep} that using RNN type architectures in~\cite{weinan2017deep} does not improve the results. Even when used in~\cite{raissi2018forward}, the RNN architecture does not improve the results, this will be shown in our work. In our new formulation, using local losses including the terminal condition helps the SGD algorithm to converge to a good local minima.

The outline of the paper is organized as follows. In the next Section, we introduce some preliminaries including the neural networks and the forward time discretization of the decoupled FBSDEs. Our deep learning-based algorithm is presented in Section~\ref{sec3}. Section~\ref{sec4} is devoted to the numerical experiments. Finally, Section~\ref{sec5} concludes this work.
\section{Preliminaries}
\label{sec2}
\subsection{The nonlinear Feynman-Kac formula}
\label{subsec21}
Let $\left(\Omega,\mathcal{F},\mathbb{P},\{\mathcal{F}_t\}_{0\le t \le T}\right)$ be a complete, filtered probability space. In this space a standard $d$-dimensional Brownian motion $W_t$ is defined, such that the filtration $\{\mathcal{F}_t\}_{0\le t\le T}$ is the natural filtration of $W_t.$ We define $|\cdot|$ as the standard Euclidean norm in the Euclidean space $\mathbb{R}$ or $\mathbb{R}^{1\times d}$ and $L^2 = L^2_{\mathcal{F}}\left(0,T; \mathbb{R}^d\right)$ the set of all $\mathcal{F}_t$-adapted and square integrable processes valued in $\mathbb{R}^d$. The triple of processes $\left(X_t,Y_t,Z_t\right):\left[0,T\right]\times\Omega\to\mathbb{R}^d \times \mathbb{R}\times\mathbb{R}^{1 \times d}$ is the solution of BSDE \eqref{eq1} if it is $\mathcal{F}_t$-adapted, square integrable, and satisfies \eqref{eq1} in the sense of
\begin{equation}
    \begin{split}
        \left\{
            \begin{array}{rcl}
                X_t & = & x_0 + \int_{0}^{t} \mu \left(s, X_s\right)\,ds + \int_{0}^{t} \sigma \left(s, X_s\right)\,dW_s,\\
   	   		    Y_t & = & g\left(X_T\right) + \int_{t}^{T} f\left(s, X_s, Y_s, Z_s\right)\,ds -\int_{t}^{T}Z_s\,dW_s,
            \end{array} \forall t \in [0, T]
        \right. 
    \end{split}
\label{eq2}
\end{equation}
where $f\left(t,X_t, Y_t,Z_t\right):\left[0,T\right]\times \mathbb{R}^d\times \mathbb{R}\times\mathbb{R}^{1\times d} \to \mathbb{R}$ is $\mathcal{F}_t$-adapted, the third term on the right-hand side is an It\^o-type integral and $g\left( X_T\right):\mathbb{R}^d \to \mathbb{R}$. This solution exist uniquely under regularity
conditions~\cite{El1997}.

One of the most important properties of BSDEs is that they provide a probabilistic representation for the solution of a specific class of partial differential equations (PDEs) given by the nonlinear Feynman–Kac formula. Consider the semi-linear parabolic PDE
\begin{equation}\label{eq3}
    \frac{\partial u}{\partial t} + \sum_{i=1}^{d} \mu_i(t, x)\frac{\partial u}{\partial x_i} + \frac{1}{2} \sum_{i, j=1}^{d} (\sigma \sigma^{\top})_{i,j}(t, x) \frac{\partial^2 u}{\partial x_{i} x_j} + f\left(t, x, u, \left(\nabla u \right) \sigma\right) = 0,
\end{equation}
with the terminal condition $u\left(T,x\right)=g(x)$. Assume that \eqref{eq3} has a classical solution $u(t, x) \in C^{1,2}([0, T] \times \mathbb{R}^d)$ and the regularity conditions of~\eqref{eq2} are satisfied. 
Then the solution of~\eqref{eq2} can be represented by 
\begin{equation}
	Y_t = u\left(t, X_t\right), \quad Z_t= \left(\nabla u\left(t, X_t\right)\right)\sigma\left(t, X_t\right) \quad \forall t \in \left[0,T\right),
\label{eq4}
\end{equation}
$\mathbb{P}$-a.s., where $\nabla u$ denotes the derivative of $u\left(t,x\right)$ with respect to the spatial variable $x$. A function approximator can be found for the solution. Due to the approximation capability in high dimensions, neural networks are a promising candidate.

\subsection{Neural Networks as function approximators}\label{subsec22}
Deep neural networks rely on the composition of simple functions, but provide an efficient way to approximate unknown functions. We introduce briefly feedforward neural networks which we will use. Let $d_0, d_1\in \mathbb{N}$ be the input and output dimensions, respectively. We fix the global number of layers as $L+2$, $L \in \mathbb{N}$ the number of hidden layers each with $n \in \mathbb{N}$ neurons. The first layer is the input layer with $d_0$ neurons and the last layer is the output layer with $d_1$ neurons. A feedforward neural network is a function $\psi_{d_0,d_1}^{\varrho,n,L}(x; \theta): \mathbb{R}^{d_0} \to \mathbb{R}^{d_1}$ as the composition
\begin{equation}\label{eq5}
    x \in \mathbb{R}^{d_0} \longmapsto T_{L+1}(\cdot;\theta^{L+1}) \circ \varrho \circ T_{L}(\cdot;\theta^{L}) \circ \varrho \circ \cdots \circ \varrho \circ T_1(x;\theta^1) \in \mathbb{R}^{d_1},
\end{equation}
where $\theta:=\left( \theta^1, \cdots, \theta^{L+1} \right) \in \mathbb{R}^{\rho}$ and $\rho$ is the number of network parameters, $x \in \mathbb{R}^{d_0}$ is the input vector. Moreover, $T_l(\cdot; \theta^l), l = 1, 2, \cdots, L+1$ are affine transformations: $T_1(x;\theta^1): \mathbb{R}^{d_0} \to \mathbb{R}^{n}$, $T_l(\cdot;\theta^l),  l = 2, \cdots, L: \mathbb{R}^{n} \to \mathbb{R}^{n}$ and $T_{L+1}(\cdot;\theta^{L+1}): \mathbb{R}^{n} \to \mathbb{R}^{d_1}$, represented by
\begin{equation*}
    T_l(y;\theta^l) = \mathcal{W}_l y + b_l,
\end{equation*}
where $\mathcal{W}_l \in \mathbb{R}^{n_{l} \times n_{l-1}}$ is the weight matrix and $b_l \in \mathbb{R}^{n_{l}}$ is the bias vector with $n_0 = d_0, n_{L+1} = d_1, n_l = n$ for $l = 1, \cdots, L$ and $\varrho: \mathbb{R} \to \mathbb{R}$ is a nonlinear function (called the activation function), and applied componentwise on the outputs of $T_l(\cdot;\theta^l)$. Common choices are $\tanh(x), \sin(x), \max\{0,x\}$ etc. The activation function must be differentiable in order to have a differentiable neural network. All the network parameters in~\eqref{eq5} given as $\theta \in \mathbb{R}^{\rho}$ can be collected as 
$$\rho = \sum_{l=1}^{L+1}n_{l}(n_{l-1}+1) = n(d_0+1) + n(n+1)(L-1) + d_1(n+1),$$
for fixed $d_0, d_1, L$ and $n$. We denote by $\Theta = \mathbb{R}^\rho$ the set of possible parameters for the neural network $\psi_{d_0,d_1}^{\varrho,n,L}(x; \theta)$ with $\theta \in \Theta$. The Universal Approximation Theorem~\cite{hornik1989multilayer} justifies the use of neural networks as function approximators.

\subsection{Learning long-term dependencies in recurrent neural networks}
\label{subsec23}
Recurrent neural networks (RNNs) are a type of artificial neural networks that allow previous outputs to be used as inputs with hidden states. It is naturally interesting to see whether RNNs can improve deep learning-based algorithms for solving BSDEs, in particular to overcome the demerits mentioned: stuck in a poor local minima or even diverge; $(Y_t, Z_t), 0<t<T$ not well approximated. However, at first glance, some advanced RNNs, e.g., Long Short-Term Memory (LSTM) networks or bidirectional RNNs should be excluded, because they do violate the markovian property for the BSDEs.

We consider the standard RNNs~\cite{rumelhart1986learning} defined as follows: given a sequence of inputs $x_1, x_2, \cdots, x_N$, each in $\mathbb{R}^{d_0}$, the network computes a sequence of hidden states $h_1, h_2, \cdots, h_N,$ each in $\mathbb{R}^n$, and a sequence of predictions $y_1, y_2, \cdots, y_N,$ each in $\mathbb{R}^{d_1}$, by the equations
\begin{equation*}
    \begin{split}
        h_i &= \varrho (\mathcal{W}_{h} h_{i-1}  + \mathcal{W}_{x} x_i + b_{h}), \\
        y_i &= \mathcal{W}_{y} h_{i} + b_{y},
    \end{split}
\end{equation*}
where $\theta:=(\mathcal{W}_{h}, \mathcal{W}_{x}, \mathcal{W}_{h}, \mathcal{W}_y, b_y) \in \mathbb{R}^\rho$ are the trainable parameters and $\varrho$ is the nonlinear activation function. Note that the standard RNNs are universal approximators as well, see~\cite{schafer2006recurrent}. If one shall think that $h_i$ depends only on the current input $x_i$ and the last hidden state $h_{i-1},$ and suppose that the distribution over the hidden states is well-defined, the standard RNNs should preserve the markovian property. However, our numerical results show that a tiny improvement can be observed.

\subsection{Forward time discretization of BSDEs}
\label{subsec24}
In order to formulate BSDEs as a learning problem, we firstly discretize the time integrals. 

The integral form of the forward SDE in~\eqref{eq1} reads
\begin{equation*}
     X_t = X_{0} + \int_0^t \mu\left(s, X_s\right)\,ds 
	          +\int_0^t \sigma\left(s, X_s\right)\,dW_s, \quad t\in \left[0,T\right].
\end{equation*}
The drift $\mu(\cdot)$ and diffusion $\sigma(\cdot)$ are assumed to be sufficiently smooth. We consider the time discretization
$$ \Delta = \{t_i|t_i \in [0, T], i = 0, 1, \cdots, N, t_i < t_{i+1}, \Delta t = t_{i+1} - t_{i}, t_0 = 0, t_N = T\}$$
for the time interval $[0,T].$ For notational convenience we write $X_i = X_{t_i}$, $W_i = W_{t_i}$, $\Delta W_i = W_{i+1} - W_i$, and $\mathcal{X}_i = \mathcal{X}_{t_i}$  for the approximations. The well-known Euler scheme reads
\begin{equation*}
     \mathcal{X}_{i+1} = \mathcal{X}_i + \mu\left(t_i, \mathcal{X}_i\right) \Delta t + \sigma\left(t_i, \mathcal{X}_i\right) \Delta W_i, \quad \text{for}\,\, i = 0, 1, \cdots, N-1,
\end{equation*}
where $\mathcal{X}_0 = X_0$ and $\Delta W_i \sim \mathcal{N}(0,\,\Delta t)$. For sufficiently small $\Delta t$, the Euler scheme has strong convergence order $\frac{1}{2}$~\cite{kloeden2013numerical}, i.e.,
\begin{equation*}
     \mathbb{E}\bigl[ | X_T - \mathcal{X}_T | \bigr] \leq C \left(\Delta t \right)^{\frac{1}{2}},
\end{equation*}
where $C > 0$ is a constant.

Next we apply the Euler scheme for the backward process. For the time interval $[t_i, t_{i+1}]$, the integral form of the backward process reads
\begin{equation*}
     Y_{t_i} = Y_{t_{i+1}} + \int_{t_{i}}^{t_{i+1}} f\left(s, X_s, Y_s, Z_s\right)\,ds -\int_{t_{i}}^{t_{i+1}} Z_s\,dW_s,
\end{equation*}
which can be straightforwardly reformulated as 
\begin{equation*}
     Y_{t_{i+1}} = Y_{t_{i}} - \int_{t_{i}}^{t_{i+1}} f\left(s, X_s, Y_s,Z_s\right)\,ds +\int_{t_{i}}^{t_{i+1}} Z_s\,dW_s.
\end{equation*}
Applying the Euler scheme for the latter equation one obtains
\begin{equation}
     \mathcal{Y}_{i+1} = \mathcal{Y}_i - f\left(t_i, \mathcal{X}_i, \mathcal{Y}_i, \mathcal{Z}_i\right) \Delta t +  \mathcal{Z}_i \Delta W_i,\quad i = 0, 1, \cdots, N-1,
\label{eq6}
\end{equation}
where $\mathcal{Y}_i = \mathcal{Y}_{t_i}$ and $\mathcal{Z}_i = \mathcal{Z}_{t_i}$ are the approximations of $Y_{t_i}$ and $Z_{t_i}$. By iterating \eqref{eq6} together with the terminal condition $g(\mathcal{X}_N)$, we have
\begin{equation}
     \mathcal{Y}_{i} = g(\mathcal{X}_N) + \sum_{j=i}^{N-1} \left(f\left(t_j, \mathcal{X}_j, \mathcal{Y}_j, \mathcal{Z}_j\right) \Delta t -  \mathcal{Z}_j \Delta W_j\right),\quad i = 0, 1, \cdots, N-1,
\label{eq7}
\end{equation}
which represents a iterative time discretization of 
\begin{equation*}
     Y_{t_i} = g(X_{T}) + \int_{t_{i}}^{T} f\left(s, X_s, Y_s, Z_s\right)\,ds -\int_{t_{i}}^{T} Z_s\,dW_s.
\end{equation*}
Note that this discretization is also used in~\cite{germain2022approximation}, their formulation is based on backward recursive local optimizations defined from~\eqref{eq7} to estimate the solution and its gradient at each time step. In our case, we consider a global optimization based on local losses obtained from~\eqref{eq7} in a forward manner. Note that the schemes in~\cite{weinan2017deep,raissi2018forward} also represent a global optimization by considering~\eqref{eq6} in a forward manner.

\section{The forward deep learning-based schemes for BSDEs}
\label{sec3}
In this section we review firstly the proposed methods in \cite{weinan2017deep, raissi2018forward}, and then present our new method.

\subsection{The deep BSDE scheme \cite{weinan2017deep}}
\label{subsec31}
The numerical approximation of $Y_{i}, i = 0, 1, \cdots, N$ in \cite{weinan2017deep} (we refer as DBSDE scheme in the rest of the paper) is designed as follows: starting from an initialization $\mathcal{Y}_0^{\theta}$ of $Y_{0}$ and $\mathcal{Z}_0^{\theta}$ of $Z_{0}$, and then using at each time step $t_i, i = 1, 2 , \cdots, N-1$ a different feedforward multilayer neural network $\psi_{d_0,d_1}^{\varrho,n,L}(x; \theta_i): \mathbb{R}^{d_0} \to \mathbb{R}^{d_1}$ to approximate $Z_{i} \in \mathbb{R}^{1\times d}$ as $\mathcal{Z}_i^{\theta}$, where the input $x$ of the network is the markovian process $\mathcal{X}_i \in \mathbb{R}^d$, $d_0 = d, d_1 = 1 \times d$. The approximation $\mathcal{Y}_i^{\theta}, i = 1, 2, \cdots, N$ is calculated using the Euler method~\eqref{eq6}. Note that this algorithm forms a global deep neural network composed of neural networks at each time step using as input data the paths of $(\mathcal{X}_i)_{i=0,1,\cdots,N}$ and $(W_i)_{i=0,1,\cdots,N}$, and gives as a final output $\mathcal{Y}_N^{\theta}$, which depends on parameters $\theta := (\mathcal{Y}_0^{\theta}, \mathcal{Z}_0^{\theta}, \theta_1, \cdots, \theta_{N-1})$. The output aims to match the terminal condition $g(\mathcal{X}_N)$ of the BSDE, and then optimizes over the parameters $\theta$ the expected square loss function:
\begin{equation*}
    \begin{split}
         \mathbf{L}(\theta) &= \mathbb{E}\bigl[|g(\mathcal{X}_N)-\mathcal{Y}_N^{\theta}|^2\bigr],\\
         \theta^{*} &\in \argmin_{\theta \in \mathbb{R}^{\rho}} \mathbf{L}(\theta),
    \end{split}
\end{equation*}
which can be done by using SGD-type algorithms. For the algorithmic framework we refer to \cite{weinan2017deep}. The DBSDE scheme uses the Adam optimizer~\cite{kingma2014adam} as an SGD optimization method with mini-batches. In the implementations, $N-1$ fully-connected feedforward neural networks are employed to approximate $\mathcal{Z}_{i}^\theta, i = 1, 2, \cdots, N-1, \theta \in \mathbb{R}^\rho$. Each of the neural networks has $L = 2$ hidden layers and $n = d+10$ neurons per hidden layer. The authors also adopt batch normalization~\cite{ioffe2015batch} right after each matrix multiplication and before activation. The rectifier function $\mathbb{R} \ni x \to \max\{0, x\} \in [0, \infty)$ is used as the activation function $\varrho$ for the hidden variables. All the weights are initialized using a normal or a uniform distribution without any pre-training. The choice of the dimension of the parameters is given as~\cite{weinan2017deep}
\begin{equation*}
    \begin{split}
        \rho &= d+1 + (N-1)(2d(d+10)+(d+10)^2+4(d+10) +2d).
    \end{split}
\end{equation*}

\subsection{The local deep BSDE scheme \cite{raissi2018forward}}
\label{subsec32}
As mentioned before, a strong drawback of the DBSDE scheme is that only $(Y_0, Z_0)$ can be well approximated. For this, \cite{raissi2018forward}~proposed to formulate the BSDE problem based on a global optimization with local losses (we refer as Local Deep BSDE or LDBSDE scheme in the rest of the paper). More precisely, the solution is approximated using a deep neural network and its gradient via AD. These approximations are performed by the global minimization of local loss functions defined from the dynamics of the BSDE at each time step given by the Euler method \eqref{eq6} and the terminal condition included as an additional term. The algorithm is given as follows:
\begin{itemize}
    \item At each time $t_i$, $i = 0, 1, 2, \cdots, N$: use one deep neural network $\psi_{d_0,d_1}^{\varrho,n,L}(x; \theta): \mathbb{R}^{d_0} \to \mathbb{R}^{d_1}$ to approximate $Y_{i} \in \mathbb{R}$ as $\mathcal{Y}_i^{\theta}$, where the input $x$ of the network is the time value $t_i \in \mathbb{R}_{+}$ and the markovian process $\mathcal{X}_i \in \mathbb{R}^d,$ $d_0 = d+1, d_1 = 1$, and
    $$\mathcal{Z}_i^{\theta} = \frac{\partial \psi_{d_0,d_1}^{\varrho,n,L}(x; \theta)}{\partial X}\Bigr|_{X = \mathcal{X}_i} \sigma \left( t_i, \mathcal{X}_i \right),$$
    a formulation based on \eqref{eq4}.
    \item The empirical loss and optimal parameters $\theta$ are given as
        \begin{equation*}
        \begin{split}
            \mathbf{L}(\theta) &= \sum_{m=1}^{M}\left( \sum_{i=0}^{N-1} |\mathcal{Y}_i^{m,\theta} - f\left(t_i, \mathcal{X}_i^m, \mathcal{Y}_i^{m,\theta}, \mathcal{Z}_i^{m,\theta}\right) \Delta t +  \mathcal{Z}_i^{m,\theta} \Delta W_i^m-\mathcal{Y}_{i+1}^{m,\theta}|^2 \right. \\
            & \quad  +\left. |\mathcal{Y}_{N}^{m,\theta} - g(\mathcal{X}_{N}^m)|^2 \vphantom{\sum_{i=0}^{N-1}} \right),\\
             \theta^{*} &\in \argmin_{\theta \in \mathbb{R}^{\rho}} \mathbf{L}(\theta),
        \end{split}
    \end{equation*}
    when using $M$ samples.
\end{itemize}
In \cite{raissi2018forward}, the author used the Adam optimizer with mini-batches, $L = 4$ hidden layers and $n = 256$ neurons. Based on this setting, the choice of the dimension of the parameters (including bias term) is given by
\begin{equation}
    \rho =  256d+198145.
\label{eq8}
\end{equation}
Furthermore, $\mathbb{R} \ni x \to \sin(x) \in [-1, 1]$ is used as activation function $\varrho$ in \cite{raissi2018forward} and the following learning rate decay approach:
$$\gamma_k = 10^{\left(\mathds{1}_{\left[20000\right]}(k) + \mathds{1}_{\left[50000\right]}(k)+\mathds{1}_{\left[80000\right]}(k) - 6\right)},$$
for $k = 1, 2, \cdots, 100000$, where $k$ is the number of the Adam optimizer steps.

\subsection{The locally additive deep BSDE scheme}
\label{subsec33}
\begin{figure}[htb!]
	\centering
	\includegraphics[width=0.7\linewidth]{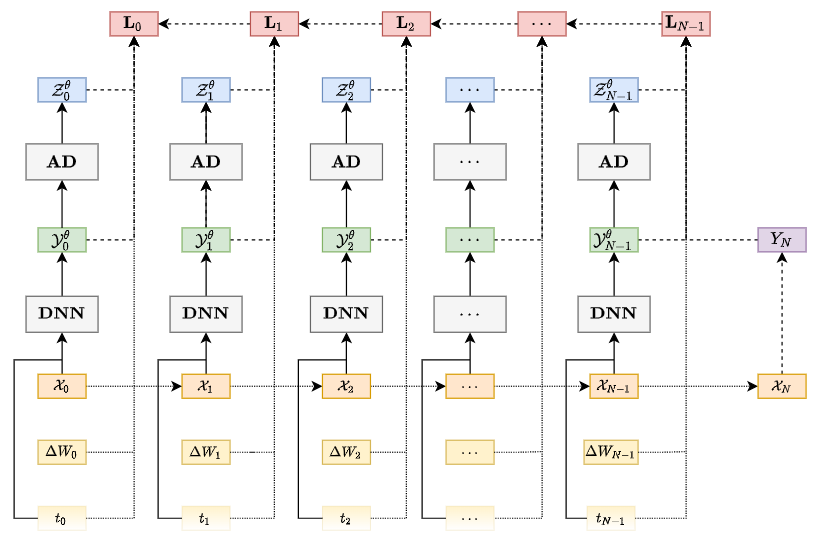}
	\caption{The architecture of the LaDBSDE scheme.}
	\label{fig1}
\end{figure}

The LDBSDE scheme improves the results of the DBSDE scheme for the approximations in the entire time domain. However, it can also get stuck in poor local minima as the DBSDE scheme especially for a complex solution structure and a long terminal time. Our idea is to consider a formulation based on a global optimization with local loss function, where each loss term includes the terminal condition. This is achieved by using the iterative time discretization \eqref{eq7}. We refer to this as the Locally additive Deep BSDE (LaDBSDE) scheme as each local loss term accumulates the information up to the terminal condition. The algorithm is given as follows:
\begin{itemize}
    \item At each time $t_i$, $i = 0, 1, 2, \cdots, N-1$: use one deep neural network $\psi_{d_0,d_1}^{\varrho,n,L}(x; \theta): \mathbb{R}^{d_0} \to \mathbb{R}^{d_1}$ to approximate $Y_{i} \in \mathbb{R}$ as $\mathcal{Y}_i^{\theta}$, where the input $x$ of the network is the time value $t_i \in \mathbb{R}_{+}$ and the markovian process $\mathcal{X}_i \in \mathbb{R}^d,$ $d_0 = d+1, d_1 = 1$, and
    $$\mathcal{Z}_i^{\theta} = \frac{\partial \psi_{d_0,d_1}^{\varrho,n,L}(x; \theta)}{\partial X}\Bigr|_{X = \mathcal{X}_i} \sigma \left( t_i, \mathcal{X}_i \right),$$
    a formulation based on~\eqref{eq4}.
    \item The empirical loss and optimal parameters $\theta$ are given as
        \begin{equation}
        \begin{split}
            \mathbf{L}(\theta) &= \sum_{m=1}^{M}\left( \sum_{i=0}^{N-1} |\mathcal{Y}_i^{m,\theta} - \sum_{j = i}^{N-1} \left(f\left(t_j, \mathcal{X}_j^m, \mathcal{Y}_j^{m,\theta}, \mathcal{Z}_j^{m,\theta}\right) \Delta t -  \mathcal{Z}_j^{m,\theta} \Delta W_j^m\right)-g(\mathcal{X}_{N}^m)|^2 \right),\\
             \theta^{*} &\in \argmin_{\theta \in \mathbb{R}^{\rho}} \mathbf{L}(\theta),
        \end{split}
        \label{eq9}
    \end{equation}
    when using $M$ samples.
\end{itemize}
We see that a neural network is used to approximate the solution of the BSDE and its gradient via AD. These approximations are obtained by the global minimization of quadratic local loss functions including the terminal time~\eqref{eq9}. The algorithmic framework (without using mini-batches and Adam optimizer) can be formulated as follows.
\begin{framework}
Let $T, \gamma \in (0, \infty)$, $d, \rho, N \in \mathbb{N}$, $X_0 \in \mathbb{R}^d$, $\mu:\left[0,T\right]\times \mathbb{R}^d \to \mathbb{R}^d$, $\sigma:\left[0,T\right]\times \mathbb{R}^d \to \mathbb{R}^{d\times d}$, $f:\left[0,T\right]\times \mathbb{R}^d\times \mathbb{R}\times\mathbb{R}^{1 \times d} \to \mathbb{R}$ and $g:\mathbb{R}^{d} \to \mathbb{R}$ be functions, let $\left(\Omega,\mathcal{F},\mathbb{P}\right)$ be a probability space, let $W^m: [0, T]\times \Omega \to \mathbb{R}^{d}$, $m \in \mathbb{N}_0$, be independent d-dimensional standard Brownian motions on $\left(\Omega,\mathcal{F},\mathbb{P}\right)$, let $t_0, t_1, \cdots, t_N \in [0, T]$ be real numbers with
\begin{equation*}
    0 = t_0 < t_1 < \cdots < t_N = T,
\end{equation*}
for every $m \in \mathbb{N}_0$ let $\mathcal{X}^m: \{ 0, 1 , \cdots, N\}\times \Omega \to \mathbb{R}^{d}$ be a stochastic process which satisfies for $i \in \{ 0, 1, \cdots, N-1\}$, $\Delta W_i^m = W_{i+1}^m - W_i^m$ that
\begin{equation*}
    \mathcal{X}_{i+1}^m = \mathcal{X}_i^m + \mu\left(t_i, \mathcal{X}_i^m\right) \Delta t + \sigma\left(t_i, \mathcal{X}_i^m\right) \Delta W_i^m, \quad \mathcal{X}_0^m = X_0,
\end{equation*}
for every $\theta \in \mathbb{R}^\rho$, $i \in \{ 0, 1 , \cdots, N-1\}$, $d_0 = d+1, d_1 = 1$, $\varrho: \mathbb{R} \to \mathbb{R}$, $L \in \mathbb{N}$, $n \in \mathbb{N}$ let $\psi_{d_0,d_1}^{\varrho,n,L}: \mathbb{R}^{d_0} \to \mathbb{R}^{d_1}$ ($\psi_{d_0,d_1}^{\varrho,n,L} \in C^1$) be a function (neural network), the output given as $\mathcal{Y}_{i}^{\theta}$ and let $\mathcal{Z}_{i}^{\theta}=\nabla_x \psi_{d_0,d_1}^{\varrho,n,L}((t_i, x);\theta)\Bigr|_{x = \mathcal{X}_i^m}\sigma(t_i, \mathcal{X}_i^m)$, for every $m \in \mathbb{N}_0$, $i \in \{ 0, 1, \cdots, N-1\}$ let $\phi^m_{i}: \mathbb{R}^\rho \times \Omega \to \mathbb{R}$ be the function which satisfies for all $\theta \in \mathbb{R}^\rho$, $\omega \in \Omega$ that
\begin{equation*}
    \phi_{i}^m(\theta, \omega) = |\mathcal{Y}_i^{m,\theta}(\omega) - \sum_{j = i}^{N-1} \left(f\left(t_j, \mathcal{X}_j^m(\omega), \mathcal{Y}_j^{m,\theta}(\omega), \mathcal{Z}_j^{m,\theta}(\omega)\right) \Delta t -  \mathcal{Z}_j^{m,\theta}(\omega) \Delta W_j^m(\omega)\right)-g(\mathcal{X}_{N}^m(\omega))|^2,
\end{equation*}
for every $m \in \mathbb{N}_0$ let $\phi^m: \mathbb{R}^\rho \times \Omega \to \mathbb{R}$ be the function which satisfies for all $\theta \in \mathbb{R}^\rho$, $\omega \in \Omega$ that
\begin{equation*}
    \phi^m(\theta, \omega) = \sum_{i=0}^{N-1} \phi_{i}^m(\theta, \omega),
\end{equation*}
for every $m \in \mathbb{N}_0$ let $\Phi^m: \mathbb{R}^\rho \times \Omega \to \mathbb{R}^\rho$ be a function which satisfies for all $\omega \in \Omega$, $\theta \in \{ v \in \mathbb{R}^\rho: ( \mathbb{R}^\rho \ni w \to \phi^m(w,\omega) \in \mathbb{R} \,\,\text{is differentiable at}\,\, v \in \mathbb{R}^\rho ) \}$ that
\begin{equation*}
    \Phi^m(\theta, \omega) = (\nabla_\theta \phi^m)(\theta, \omega),
\end{equation*}
and let $\varTheta: \mathbb{N}_0 \times \Omega \to \mathbb{R}^\rho$ be a stochastic process which satisfy for all $m \in \mathbb{N}$ that
\begin{equation*}
    \varTheta_m = \varTheta_{m-1} - \gamma \Phi^m(\varTheta_{m-1}).
\end{equation*}
\label{frame1}
\end{framework}
The architecture of the LaDBSDE scheme is displayed in Figure~\ref{fig1}. 
The flow of the information is represented by the direction of the arrows. The calculations can be broken down into three steps. In the first step, the samples of the forward SDE are calculated. The information used in this step is represented by the dotted lines. For instance, to calculate $\mathcal{X}_2$, $(t_1, \Delta W_1, \mathcal{X}_1)$ is used, and $(t_{N-1}, \Delta W_{N-1}, \mathcal{X}_{N-1})$ for $\mathcal{X}_N$. The second step is to calculate the values $(\mathcal{Y}_i^{\theta}, \mathcal{Z}_i^{\theta})$ for $i = 0, 1, \cdots, N-1$ using a deep neural network (DNN) and the AD. The information needed for such calculations is represented by the solid lines. For example, the DNN uses as input $(t_1, \mathcal{X}_1)$ to calculate $\mathcal{Y}_1^{\theta}$. Using the AD we calculate the gradient in the spatial direction to obtain $\mathcal{Z}_1^{\theta}$. Finally, the local losses are calculated backwardly with the information presented by the dashed lines. To calculate $\mathbf{L}_{N-1}$, the terminal condition $Y_N = g(\mathcal{X}_N)$ and $(t_{N-1}, \Delta W_{N-1}, \mathcal{X}_{N-1}, \mathcal{Y}_{N-1}^{\theta}, \mathcal{Z}_{N-1}^{\theta})$ are used. For $\mathbf{L}_{N-2}$, $(t_{N-2}, \Delta W_{N-2}, \mathcal{X}_{N-2}, \mathcal{Y}_{N-2}^{\theta}, \mathcal{Z}_{N-2}^{\theta})$ and the information form $\mathbf{L}_{N-1}$ are used, namely $Y_N$ and $(t_{N-1}, \Delta W_{N-1}, \mathcal{X}_{N-1}, \mathcal{Y}_{N-1}^{\theta})$. The same holds for the other loss terms. We use a backward implementation of the local loss functions because it is more efficient than their forward implementation. The forward and backward implementations of~\eqref{eq9} for one sample are given in Algorithm~\ref{alg1} and~\ref{alg2}, respectively.
\begin{algorithm}
\begin{algorithmic}
\caption{A forward implementation of the loss function~\eqref{eq9}}\label{alg1}
\State \textbf{Data}: $(t_i, \Delta W_i, \mathcal{X}_i, \mathcal{Y}^{\theta}_i, \mathcal{Z}^{\theta}_i)_{0\leq i \leq N-1}$, $\mathcal{X}_N$, $\Delta t$
\State \textbf{Result:} $\mathbf{L}$
\State $\mathbf{L} \gets 0$
\For{$i= 0: N-1$}
    \State $\tilde{\mathcal{Y}}^{\theta}_i \gets g(\mathcal{X}_N)$
    \For{$j= i: N-1$}
        \State $\tilde{\mathcal{Y}}^{\theta}_i \gets \tilde{\mathcal{Y}}^{\theta}_i + f(t_j, \mathcal{X}_j, \mathcal{Y}^{\theta}_j, \mathcal{Z}^{\theta}_j) \Delta t - \mathcal{Z}^{\theta}_j \Delta W_j$
    \EndFor
    \State $\mathbf{L}_i \gets |\mathcal{Y}^{\theta}_i - \tilde{\mathcal{Y}}^{\theta}_i|^2$
    \State $\mathbf{L} \gets \mathbf{L} + \mathbf{L}_i$
\EndFor
\end{algorithmic}
\end{algorithm}
\begin{algorithm}
\caption{A backward implementation of the loss function~\eqref{eq9}}\label{alg2}
\begin{algorithmic}
\State \textbf{Data}: $(t_i, \Delta W_i, \mathcal{X}_i, \mathcal{Y}^{\theta}_i, \mathcal{Z}^{\theta}_i)_{0\leq i \leq N-1}$, $\mathcal{X}_N$, $\Delta t$
\State \textbf{Result:} $\mathbf{L}$
\State $\mathbf{L} \gets 0$
\State $\tilde{\mathcal{Y}}^{\theta}_N \gets g(\mathcal{X}_N)$
\For{$i= N-1:0$}
    \State $\tilde{\mathcal{Y}}^{\theta}_i \gets \tilde{\mathcal{Y}}^{\theta}_{i+1} + f(t_i, \mathcal{X}_i, \mathcal{Y}^{\theta}_i, \mathcal{Z}^{\theta}_i) \Delta t - \mathcal{Z}^{\theta}_i \Delta W_i$
\EndFor
\For{$i= 0:N-1$}
    \State $\mathbf{L}_i \gets |\mathcal{Y}^{\theta}_i - \tilde{\mathcal{Y}}^{\theta}_i|^2$
    \State $\mathbf{L} \gets \mathbf{L} + \mathbf{L}_i$
\EndFor
\end{algorithmic}
\end{algorithm}
With Algorithm ~\ref{alg1} the computation time of LaDBSDE is comparable to that of LDBSDE.

We consider a similar network architecture as in~\cite{raissi2018forward}. Based on Framework~\ref{frame1} we require to optimize over differentiable deep neural networks, and using the classical rectifier function may lead to an explosion while calculating the numerical approximation of the $Z$ process. We consider $\mathbb{R} \ni x \to \tanh(x) \in [-1, 1]$. Moreover, using $L = 4$ hidden layers and $n = 10+d$ neurons for the hidden layers is enough, increasing $L$ or $n$ does not improve the accuracy in our tests. The dimension of the parameters is given in Remark~\ref{r1}.
\begin{remark}
Let $\rho \in \mathbb{N}$ be the dimension of the parameters in LaDBSDE scheme.
\begin{enumerate}
    \item $(10+d)(d+1+1)$ components of $\theta \in \mathbb{R}^\rho$ are used to uniquely describe the linear transformation from (d+1)-dimensional input layer to (10+d)-dimensional first hidden layer.
    \item $(10+d)(10+d+1)$ components of $\theta \in \mathbb{R}^\rho$ are used to uniquely describe the linear transformation from (10+d)-dimensional first hidden layer to (10+d)-dimensional second hidden layer.
    \item $(10+d)(10+d+1)$ components of $\theta \in \mathbb{R}^\rho$ are used to uniquely describe the linear transformation from (10+d)-dimensional second hidden layer to (10+d)-dimensional third hidden layer.
    \item $(10+d)(10+d+1)$ components of $\theta \in \mathbb{R}^\rho$ are used to uniquely describe the linear transformation from (10+d)-dimensional third hidden layer to (10+d)-dimensional fourth hidden layer.
    \item $10+d+1$ components of $\theta \in \mathbb{R}^\rho$ are used to uniquely describe the linear transformation from (10+d)-dimensional fourth hidden layer to 1-dimensional output layer.
\end{enumerate}
Therefore, $\rho$ is given as
\begin{equation}
    \begin{split}
        \rho &=  \underbrace{(10+d)(d+1+1)}_\text{item 1.} + \underbrace{3(10+d)(10+d+1)}_\text{items 2.-4.}+\underbrace{(10+d+1)}_\text{item 5.}\\
        &= 2d^2+56d+361.
    \end{split}
\label{eq10}
\end{equation}
\label{r1}
\end{remark}
Compared to the complexity~\eqref{eq8} given in~\cite{raissi2018forward}, our parametrization of the neural network gives a smaller complexity~\eqref{eq10}. For instance, considering an example in $d=100$, the complexity based on equation \eqref{eq10} is decreased with a factor around $9$ when compared to \eqref{eq8}. In order to further reduce the computation time compared to the learning approach given in~\cite{raissi2018forward}, we consider a learning rate decay optimization approach based on the relative magnitude of the loss function \cite{chan2019machine}. We start with a learning rate $\gamma_0$. For each $1000$ optimization steps, we evaluate the loss every $100$ steps on a validation size of $1024$. Then we can take the average of $10$ collected loss values. If the relative loss over two consecutive periods is less than $5\%$, we have reached a loss plateau and reduce the learning rate by half. To avoid using very small learning rates, we set a threshold $ \gamma_{min}$. If the loss value doesn't decrease any more, the learning process is terminated. Otherwise, we continue until $60000$ optimization steps.
The hyperparameter values used for all schemes are reported in Table~\ref{tab1}, which give the best approximations in each scheme in our numerical experiments.
\begin{table}[h!]
{\footnotesize
\begin{center}
  \begin{tabular}{| c | c | c | c | c | c | c |}
  \hline
   \multirow{2}{*}{Scheme} & \multicolumn{4}{c|}{Network parametrization} & \multicolumn{2}{c|}{Learning rate decay}\\ \cline{2-7}
    & \# networks & L & n & $\varrho$ & $\gamma_0$ & $\gamma_{min}$\\ \hline
    DBSDE & N-1 & 2 & 10+d & $\max(0, x)$ & $10^{-2}$ & $10^{-4}$\\ \hline
    LDBSDE & 1 & 4 & 10+d & $\sin(x)$ & $10^{-3}$ & $10^{-5}$\\ \hline
    LaDBSDE & 1 & 4 & 10+d & $\tanh(x)$ & $10^{-3}$ & $10^{-5}$\\ \hline
   \end{tabular}
  \end{center}
\caption{Hyperparameters for all the schemes.}
\label{tab1}  
}
\end{table}

\section{Numerical results}
\label{sec4}
In this section we illustrate the improved performance using the LaDBSDE scheme compared to the schemes DBSDE and LDBSDE. The results are presented using $10$ independent runs of the algorithms.

We start with an example where the DBSDE method diverges.

\begin{exmp}
Consider the decoupled FBSDE \cite{hure2020deep}
\begin{equation*}
    \begin{split}
        \left\{
            \begin{array}{rcl}
                dX_t & = & \mu \,dt + \sigma \, dW_t, \quad X_0 = x_0,\\ 
                -dY_t & = & \left(\left( \cos\left(\bar{X}\right)  + 0.2 \sin\left(\bar{X}\right)\right)\exp\left(\frac{T-t}{2}\right) \right. \\
                & & \left. -\frac{1}{2} \left( \sin\left(\bar{X}\right) \cos\left(\bar{X}\right) \exp\left( T-t\right) \right)^2  + \frac{1}{2}\left(Y_t \bar{Z}\right)^2 \vphantom{\cos} \right)\,dt -Z_t \,dW_t,\\  
   	   	    	Y_T & = & \cos\left(\bar{X}\right),
            \end{array}
        \right. 
    \end{split}
\end{equation*}
\label{ex1}
\end{exmp}
where $\bar{X} = \sum_{i=1}^{d}X_t^i$ and $\bar{Z} =\sum_{i=1}^{d}Z_t^i$. The analytical solution is given by
\begin{equation*}
    \begin{split}
        \left\{
            \begin{array}{rcl}
                Y_t & = & \exp\left(\frac{T-t}{2}\right) \cos\left( \bar{X} \right),\\
   	   	    	Z_t & = &-\sigma \exp\left(\frac{T-t}{2}\right) \sin\left(\bar{X} \right)\mathds{1}_{\mathbb{R}^d}.
            \end{array}
        \right. 
    \end{split}
\end{equation*}

We begin with $d=1$, the exact solution of $\left(Y_0, Z_0\right) \doteq \left( 1.4687, -2.2874\right)$ for $T=2$, $\mu = 0.2$, $\sigma = 1$ and $x_0 = 1$. A testing sample of $4096$ and $30000$ optimization steps are used. Firstly, we test the approximations of $Y_0$ and $Z_0$ of all the schemes by comparing the mean absolute errors defined as $\bar{\epsilon}_{Y_0} = \frac{1}{10}\sum_{i=1}^{10}| Y_0 - \mathcal{Y}_0^{i,\theta} |$ and $\bar{\epsilon}_{Z_0} = \frac{1}{10}\sum_{i=1}^{10} \left(\frac{1}{d} \sum_{j=1}^{d}| Z_0^j - \mathcal{Z}_0^{i,j,\theta} |\right)$. The results are reported in Table~\ref{tab2} by varying $N.$ 
\begin{table}[h!]
{\footnotesize
\begin{center}
  \begin{tabular}{| c | c | c | c | c |}
  \hline
   \multirow{3}{*}{Scheme} & N = 120 & N = 160 & N = 200 & N = 240\\ \
    &  $\bar{\epsilon}_{Y_0}$ (Std. Dev.) & $\bar{\epsilon}_{Y_0}$ (Std. Dev.) & $\bar{\epsilon}_{Y_0}$ (Std. Dev.) & $\bar{\epsilon}_{Y_0}$ (Std. Dev.)\\ 
    &  $\bar{\epsilon}_{Z_0}$ (Std. Dev.) & $\bar{\epsilon}_{Z_0}$ (Std. Dev.) & $\bar{\epsilon}_{Z_0}$ (Std. Dev.) & $\bar{\epsilon}_{Z_0}$ (Std. Dev.)\\ \hline
   \multirow{2}{*}{DBSDE} & NC & NC & NC & NC \\ 
    & NC & NC & NC & NC \\ \hline
   \multirow{2}{*}{LDBSDE} & 5.08e-1 (1.85e-1)  & 4.70e-1 (1.72e-1) & 4.55e-1 (1.69e-1) & 4.46e-1 (1.65e-1) \\
   & 2.88e-1 (1.59e-1)  & 2.35e-1 (1.54e-1) & 2.12e-1 (1.45e-1) & 2.00e-1 (1.43e-1) \\ \hline
   \multirow{2}{*}{LDBSDE (RNN)} & 6.33e-1 (2.68e-1)  & 5.08e-1 (1.99e-1) & 4.79e-1 (2.39e-1) & 4.44e-1 (2.50e-1) \\
     & 3.70e-1 (1.87e-1)  & 2.79e-1 (1.20e-1) & 3.28e-1 (1.85e-1) & 2.64e-1 (2.12e-1) \\ \hline
   \multirow{2}{*}{LDBSDE (LSTM)} & 8.85e-1 (9.90e-2)  & 8.00e-1 (9.15e-2) & 8.01e-1 (1.12e-1) & 7.28e-1 (1.09e-1) \\
     & 4.87e-1 (6.96e-2)  & 4.23e-1 (5.31e-2) & 4.40e-1 (8.38e-2) & 3.86e-1 (7.25e-2) \\ \hline
   \multirow{2}{*}{LaDBSDE} & 1.17e-1 (3.94e-2) & 1.01e-1 (3.26e-2) & 8.66e-2 (2.86e-2) & 7.90e-2 (2.69e-2)\\
   & 5.98e-2 (3.26e-2) & 5.83e-2 (3.08e-2) & 5.41e-2 (2.71e-2) & 4.94e-2 (2.86e-2)\\ \hline   
   \end{tabular}
  \end{center}
\caption{The mean absolute errors of $Y_0$ and $Z_0$ for Example~\ref{ex1} using $d = 1$. The standard deviation is given in parenthesis.}
\label{tab2}  
}
\end{table}
\begin{figure}[h!]
	\centering
	\begin{subfigure}[h!]{0.48\linewidth}
		\includegraphics[width=\linewidth]{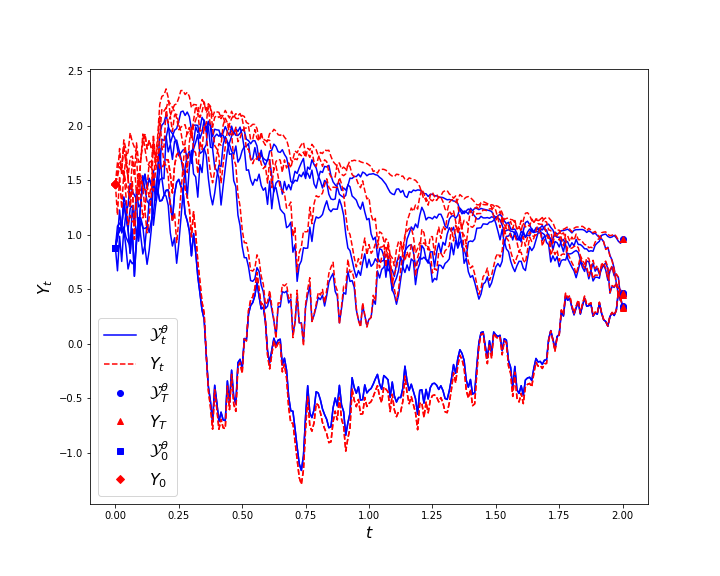}
		\caption{LDBSDE $Y$ samples.}
		\label{fig2a}
	\end{subfigure}
	\begin{subfigure}[h!]{0.48\linewidth}
		\includegraphics[width=\linewidth]{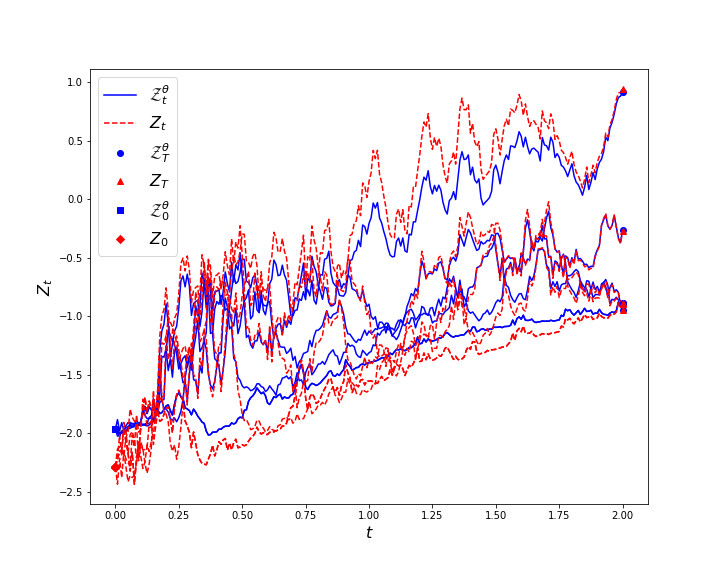}
		\caption{LDBSDE $Z$ samples.}
		\label{fig2b}
	\end{subfigure}
	\begin{subfigure}[h!]{0.48\linewidth}
		\includegraphics[width=\linewidth]{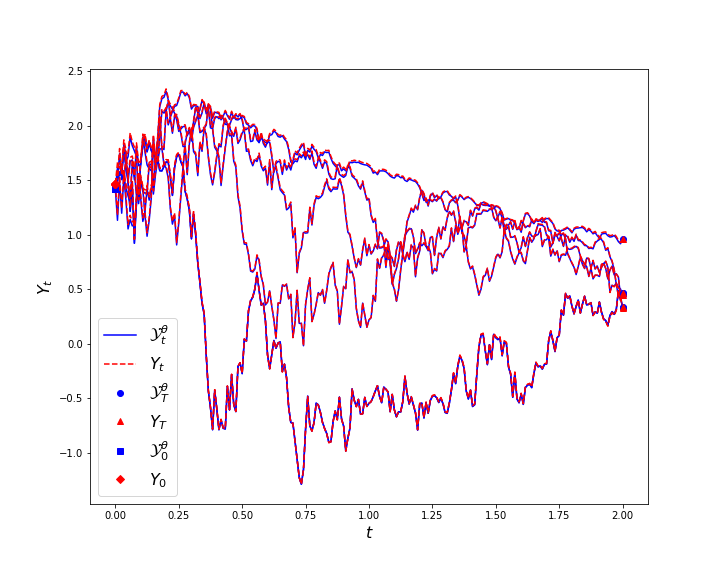}
		\caption{LaDBSDE $Y$ samples.}
		\label{fig2c}
	\end{subfigure}
	\begin{subfigure}[h!]{0.48\linewidth}
		\includegraphics[width=\linewidth]{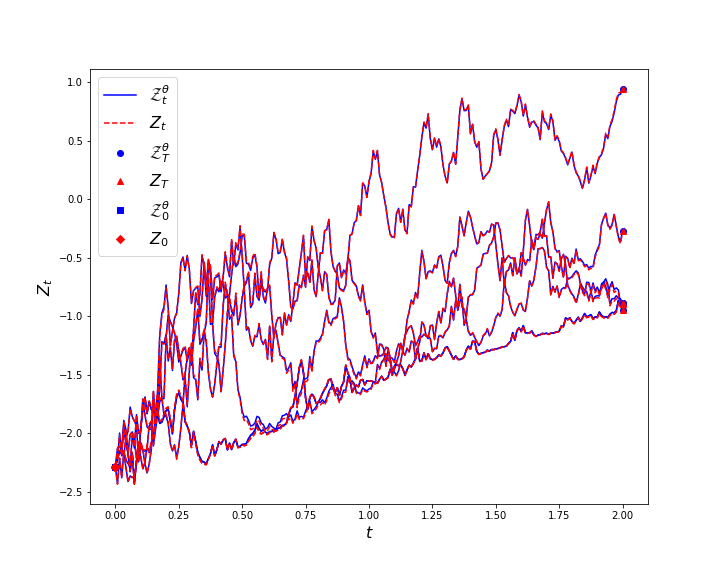}
		\caption{LaDBSDE $Z$ samples.}
		\label{fig2d}
	\end{subfigure}	
    \caption{Realizations of $5$ independent paths for Example~\ref{ex1} using $d = 1$ and $N=240$. $(Y_t, Z_t)$ and $(\mathcal{Y}_t^{\theta},\mathcal{Z}_t^{\theta})$ are exact and learned solutions for $t \in [0, T]$, respectively.}
	\label{fig2}
\end{figure}
Actually, only a few hundreds optimization steps are needed to achieve a good approximation of $(Y_0, Z_0).$ However, to obtain good approximations for the whole time domain, a high number of optimization steps is needed. From Table~\ref{tab2} we see that the DBSDE scheme diverges. The LDBSDE scheme converges to a poor local minima, the relative errors with $N = 240$ are around $30.37\%$ and $8.74\%$ for $Y_0$ and $Z_0$ respectively. In order to numerically test that the RNN type architectures does not help the LDBSDE scheme to overcome the issue of poor local minima, we use the RNN and LSTM architectures, which are referred as LDBSDE (RNN) and LDBSDE (LSTM), respectively. Using the LSTM architecture in the LDBSDE scheme, the approximation errors are quite high, since the LSTM violates the markovian property of the BSDEs. Even using the RNN in the LDBSDE scheme cannot improve the approximations. The LaDBSDE scheme gives smaller relative errors than the LDBSDE, $5.38\%$ and $2.16\%$ for $Y_0$ and $Z_0$, respectively. Note that the approximation of $Y_0$ in \cite{hure2020deep} is more accurate than all the schemes (the results for $Z_0$ are missing) in this example. However, the algorithm in~\cite{hure2020deep} is a backward type approach, which is based on local optimizations at each time step. Its computational cost should be much higher than all the DBSDE, LDBSDE and LaDBSDE schemes.

\begin{figure}[h!]
	\centering
	\begin{subfigure}[h!]{0.48\linewidth}
		\includegraphics[width=\linewidth]{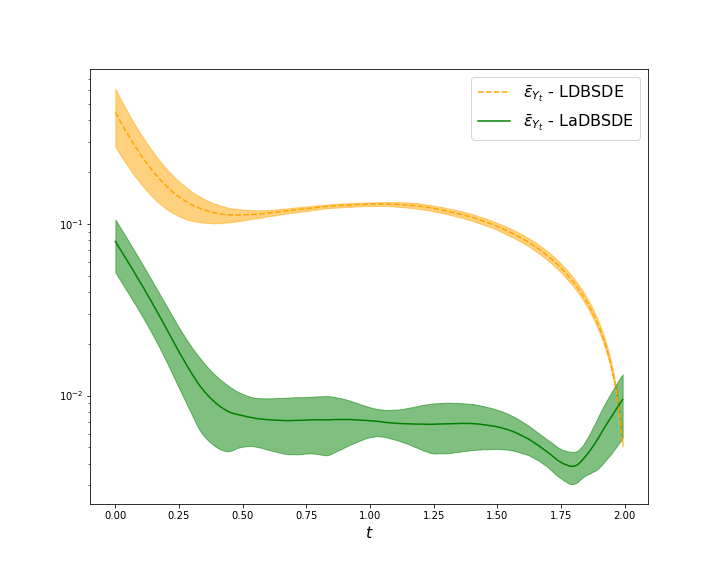}
		\caption{$Y$ process.}
		\label{fig3a}
	\end{subfigure}
	\begin{subfigure}[h!]{0.48\linewidth}
		\includegraphics[width=\linewidth]{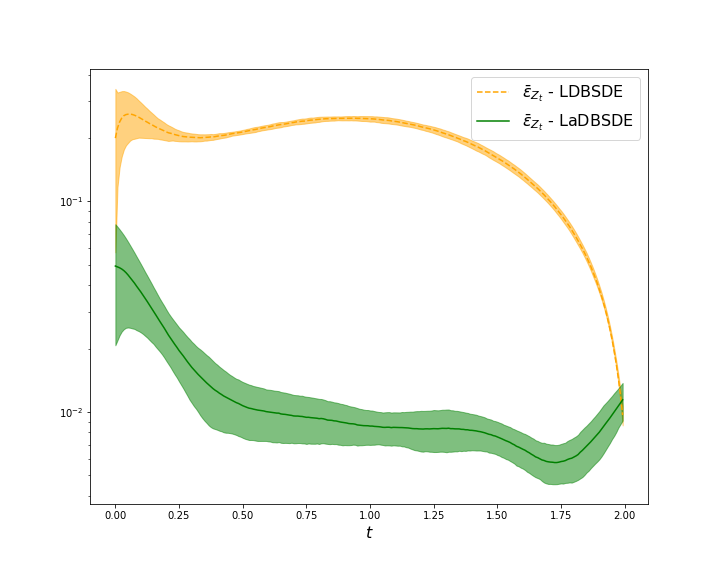}
		\caption{$Z$ process.}
		\label{fig3b}
	\end{subfigure}
	\caption{The mean regression errors  $(\bar{\epsilon}_{Y_i}, \bar{\epsilon}_{Z_i})$ at time step $t_i, i = 0, \cdots, N-1$ for Example~\ref{ex1} using $d=1$ and $N=240$. The standard deviation is given in the shaded area.}
	\label{fig3}
\end{figure}

Next we compare the performances of LDBSDE and LaDBSDE for the entire time domain. We display $5$ paths of processes $Y$ and $Z$ with $N = 240$ in Figure~\ref{fig2}. Note that the approximation for the entire time domain is not discussed in \cite{weinan2017deep}, and in \cite{raissi2018forward} only $Y$ is considered. From Figure~\ref{fig2} we see that LaDBSDE outperforms the LDBSDE scheme. In order to evaluate the accuracy at each time step for all the testing sample of $4096$, we calculate the mean regression errors defined as $\bar{\epsilon}_{Y_i} = \frac{1}{10}\sum_{j=1}^{10}\mathbb{E}[| Y_i - \mathcal{Y}_i^{j,\theta} |]$ and $\bar{\epsilon}_{Z_i} = \frac{1}{10}\sum_{j=1}^{10}\left(\frac{1}{d} \sum_{l=1}^{d}\mathbb{E}[| Z_i^l - \mathcal{Z}_i^{l,j,\theta} |]\right)$ for $i = 0, 1, \cdots, N-1.$ The results are presented in Figure~\ref{fig3}. 
We see that LaDBSDE scheme gives smaller regression errors at each time layer.

\begin{table}[h!]
{\footnotesize
\begin{center}
  \begin{tabular}{| c | c | c | c | c |}
  \hline
  \multirow{3}{*}{Scheme} & N = 60 & N = 80 & N = 100 & N = 120\\ \
    &  $\bar{\epsilon}_{Y_0}$ (Std. Dev.) & $\bar{\epsilon}_{Y_0}$ (Std. Dev.) & $\bar{\epsilon}_{Y_0}$ (Std. Dev.) & $\bar{\epsilon}_{Y_0}$ (Std. Dev.)\\ 
    &  $\bar{\epsilon}_{Z_0}$ (Std. Dev.) & $\bar{\epsilon}_{Z_0}$ (Std. Dev.) & $\bar{\epsilon}_{Z_0}$ (Std. Dev.) & $\bar{\epsilon}_{Z_0}$ (Std. Dev.)\\ \hline
    \multirow{2}{*}{DBSDE} & 5.89e-2 (1.24e-3) & 6.00e-2 (1.81e-3) & 6.09e-2 (1.92e-3) & 6.25e-2 (1.98e-3) \\ 
    & 7.46e-3 (4.04e-4) & 6.39e-3 (4.82e-4) & 5.85e-3 (5.38e-4) & 5.58e-3 (5.29e-4) \\ \hline
   \multirow{2}{*}{LDBSDE} & 8.34e-2 (1.39e-2) & 9.57e-2 (1.88e-2) & 9.61e-2 (8.94e-3) & 9.17e-2 (1.31e-2) \\ 
    & 5.91e-3 (1.30e-3) & 7.67e-3 (4.10e-3) & 6.57e-3 (1.53e-3) & 5.93e-3 (1.17e-3) \\ \hline
   \multirow{2}{*}{LaDBSDE} & 1.94e-2 (2.61e-2) & 9.14e-3 (6.35e-3) & 7.30e-3 (6.23e-3) & 6.20e-3 (4.65e-3)\\ 
   & 3.54e-3 (1.22e-3) & 2.97e-3 (5.17e-4) & 4.09e-3 (1.04e-3) & 3.11e-3 (1.28e-3)\\ \hline
   \end{tabular}
  \end{center}
\caption{The mean absolute errors of $Y_0$ and $Z_0$ for Example~\ref{ex1} using $d = 100$. The standard deviation is given in parenthesis.}
\label{tab3} 
}
\end{table}
\begin{figure}[h!]
	\centering
	\begin{subfigure}[h!]{0.48\linewidth}
		\includegraphics[width=\linewidth]{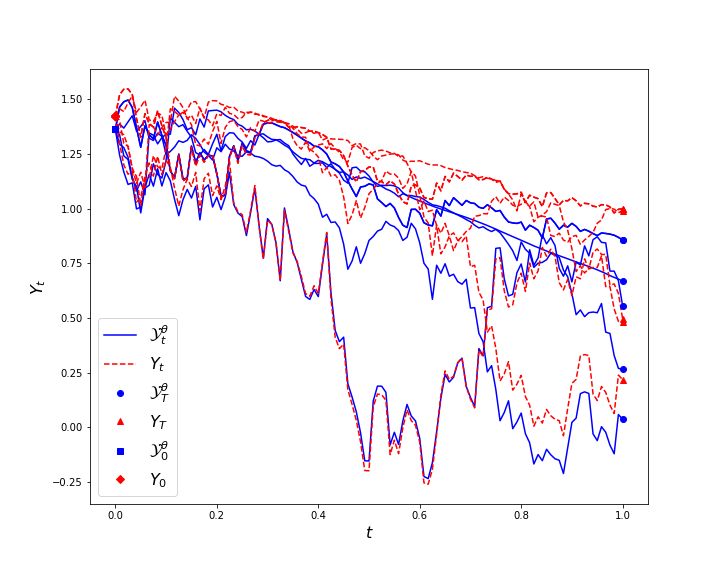}
		\caption{DBSDE $Y$ samples.}
		\label{fig4a}
	\end{subfigure}
	\begin{subfigure}[h!]{0.48\linewidth}
		\includegraphics[width=\linewidth]{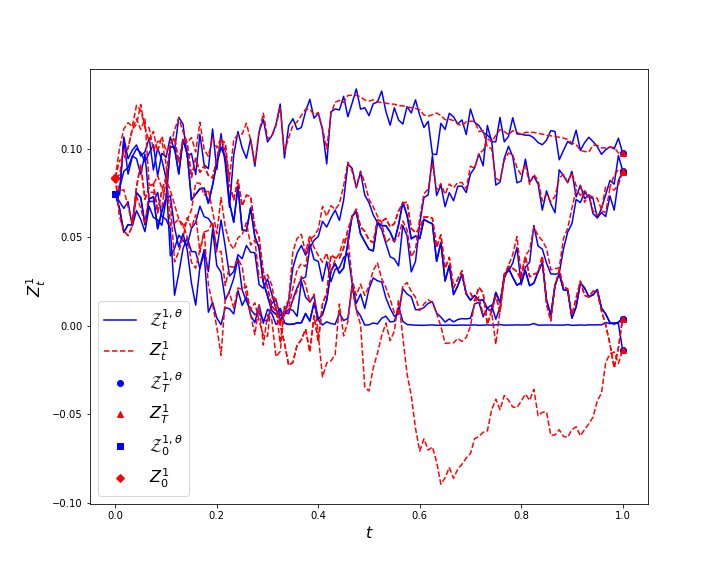}
		\caption{DBSDE $Z^1$ samples.}
		\label{fig4b}
	\end{subfigure}
	\begin{subfigure}[h!]{0.48\linewidth}
		\includegraphics[width=\linewidth]{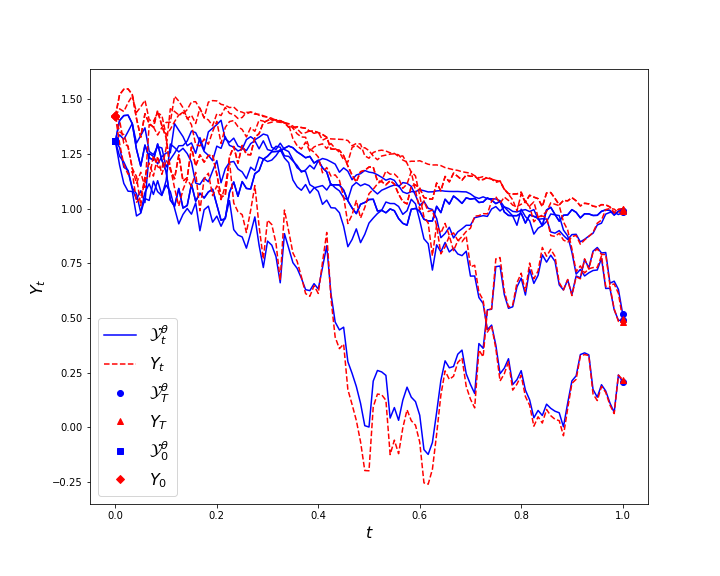}
		\caption{LDBSDE $Y$ samples.}
		\label{fig4c}
	\end{subfigure}
	\begin{subfigure}[h!]{0.48\linewidth}
		\includegraphics[width=\linewidth]{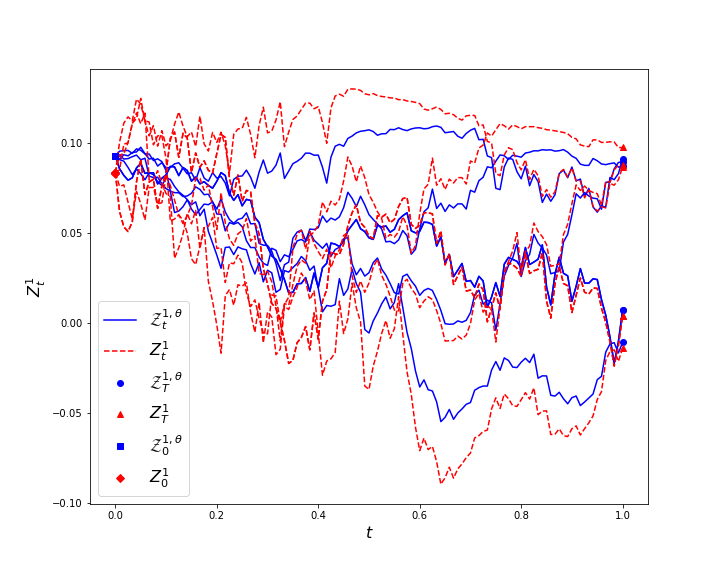}
		\caption{LDBSDE $Z^1$ samples.}
		\label{fig4d}
	\end{subfigure}	
	\begin{subfigure}[h!]{0.48\linewidth}
		\includegraphics[width=\linewidth]{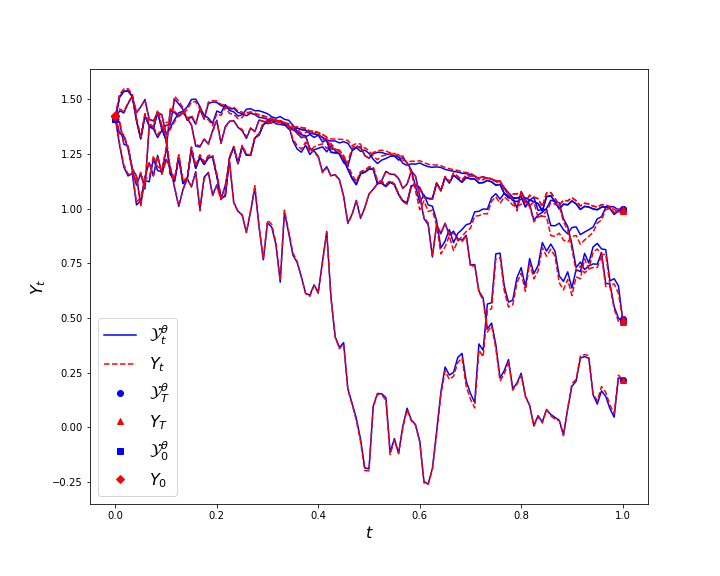}
		\caption{LaDBSDE $Y$ samples.}
		\label{fig4e}
	\end{subfigure}	
	\begin{subfigure}[h!]{0.48\linewidth}
		\includegraphics[width=\linewidth]{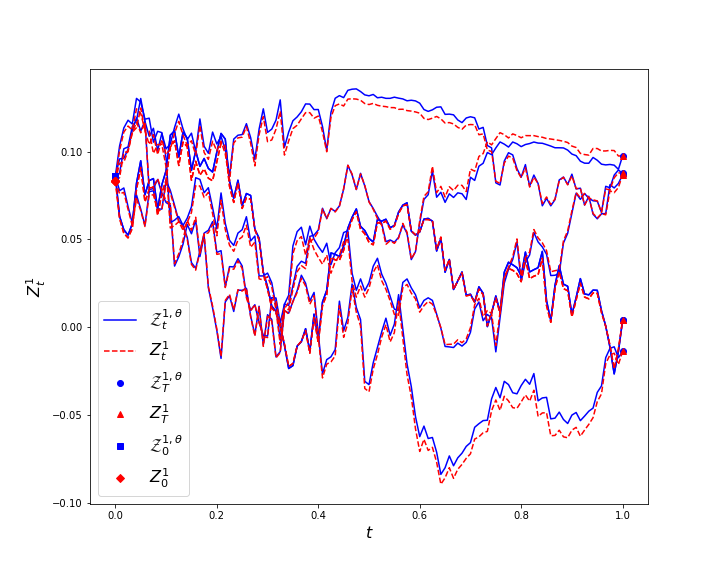}
		\caption{LaDBSDE $Z^1$ samples.}
		\label{fig4f}
	\end{subfigure}		
	 \caption{Realizations of $5$ independent paths for Example~\ref{ex1} using $d = 100$ and $N=120$. $(Y_t, Z_t^1)$ and $(\mathcal{Y}_t^{\theta},\mathcal{Z}_t^{1,\theta})$ are exact and learned solutions for $t \in [0, T]$, respectively.}
	\label{fig4}
\end{figure}
We consider the high dimensional case by setting $d=100.$ The exact solution for $T=1$, $\mu = \frac{0.2}{d}$, $\sigma = \frac{1}{\sqrt{d}}$ and $x_0 = 1$ is $\left(Y_0, Z_0\right) \doteq \left( 1.4217, (0.0835, \cdots, 0.0835)\right)$. Here we use $60000$ optimization steps. The numerical approximation of each scheme is reported in Table~\ref{tab3} for $Y_0$ and $Z_0$ by varying $N$. In contract to the one-dimensional case, we observe that the DBSDE scheme gives good approximations in this example for $d=100$ and maturity $T=1$. The reason could be that the diffusion reduces due to the large value of dimensionality ($\sigma = \frac{1}{\sqrt{d}}$), and the maturity is shorter than that in the case of one dimension. The DBSDE scheme diverges by setting $T=2$. Nevertheless, the smallest errors are still given by the LaDBSDE scheme. 

\begin{figure}[h!]
	\centering
	\begin{subfigure}[h!]{0.48\linewidth}
		\includegraphics[width=\linewidth]{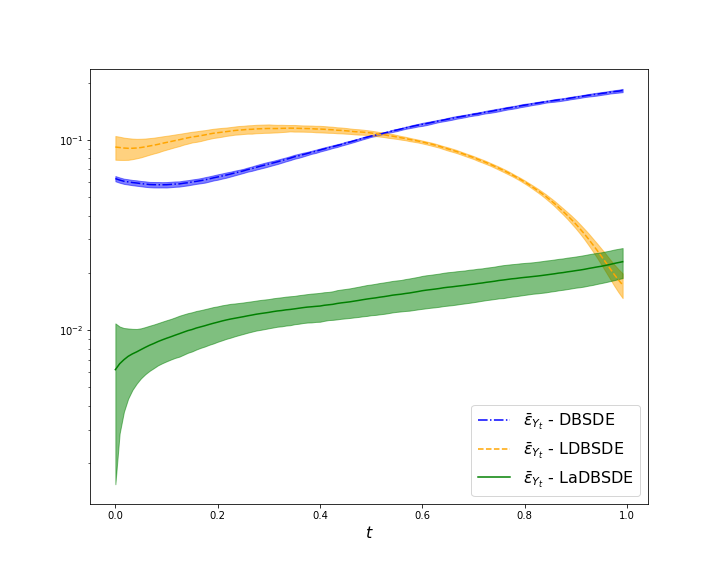}
		\caption{$Y$ process.}
		\label{fig5a}
	\end{subfigure}
	\begin{subfigure}[h!]{0.48\linewidth}
		\includegraphics[width=\linewidth]{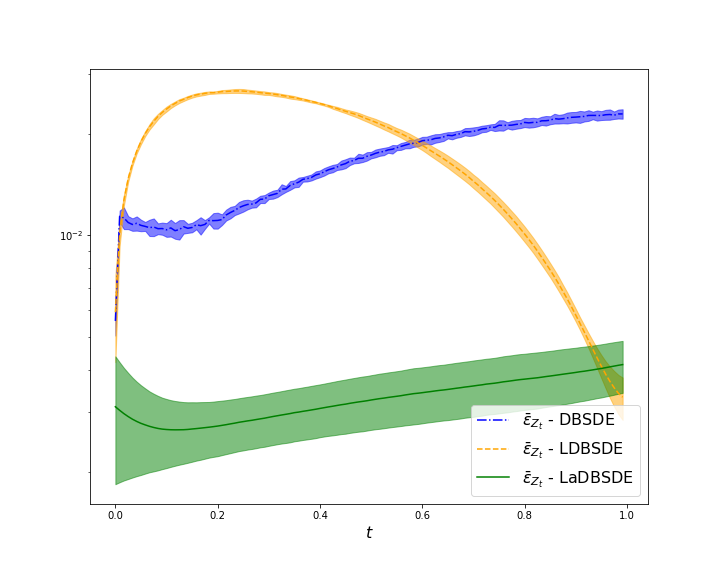}
		\caption{$Z$ process.}
		\label{fig5b}
	\end{subfigure}
	\caption{The mean regression errors  $(\bar{\epsilon}_{Y_i}, \bar{\epsilon}_{Z_i})$ at time step $t_i, i = 0, \cdots, N-1$ for Example~\ref{ex1} using $d=100$ and $N=120$. The standard deviation is given in the shaded area.}
	\label{fig5}
\end{figure}

To compare the approximations for the entire time domain in the high dimensional case we display $5$ paths with $N = 120$ of process $Y$ and the first component of $Z$ in Figure~\ref{fig4}. Note that the approximation quality of the other components in $Z$ is the same as that of $Z^1$. The regression errors are given in Figure~\ref{fig5}. 
Our method shows better approximations of processes $Y$ and $Z$ on the entire time domain compared to the DBSDE and LDBSDE schemes.

Next we consider the example with a driver function in which the $Z$ process grows quadratically.
\begin{table}[h!]
{\footnotesize
\begin{center}
  \begin{tabular}{| c | c | c | c | c |}
  \hline
  \multirow{3}{*}{Scheme} & N = 60 & N = 80 & N = 100 & N = 120\\ \
    &  $\bar{\epsilon}_{Y_0}$ (Std. Dev.) & $\bar{\epsilon}_{Y_0}$ (Std. Dev.) & $\bar{\epsilon}_{Y_0}$ (Std. Dev.) & $\bar{\epsilon}_{Y_0}$ (Std. Dev.)\\ 
    &  $\bar{\epsilon}_{Z_0}$ (Std. Dev.) & $\bar{\epsilon}_{Z_0}$ (Std. Dev.) & $\bar{\epsilon}_{Z_0}$ (Std. Dev.) & $\bar{\epsilon}_{Z_0}$ (Std. Dev.)\\ \hline    \multirow{2}{*}{DBSDE} & 8.89e-4 (3.65e-4) & 8.17e-4 (3.50e-4) & 1.09e-3 (3.73e-4) & 8.88e-4 (5.11e-4) \\
    & 8.05e-4 (5.66e-5) & 8.40e-4 (7.63e-5) & 9.55e-4 (1.06e-4) & 9.99e-4 (7.77e-5) \\ \hline
   \multirow{2}{*}{LDBSDE} & 1.45e-3 (6.64e-4) & 1.67e-3 (7.64e-4) & 3.03e-3 (2.66e-3) & 3.31e-4 (2.72e-3) \\
    & 3.35e-4 (9.22e-5) & 4.39e-4 (1.99e-4) & 5.87e-4 (2.10e-4) & 5.25e-4 (2.06e-4) \\ \hline
   \multirow{2}{*}{LaDBSDE} & 6.95e-4 (3.49e-4) & 7.62e-4 (5.59e-4) & 5.92e-4 (3.50e-4) & 9.98e-4 (5.75e-4)\\
   & 1.52e-4 (2.14e-5) & 1.58e-4 (3.56e-5) & 1.42e-4 (3.44e-5) & 1.74e-4 (4.80e-5) \\ \hline
   \end{tabular}
  \end{center}
 \caption{The mean absolute errors of $Y_0$ and $Z_0$ for Example~\ref{ex2} using $d = 100$. The standard deviation is given in parenthesis.}
\label{tab4} 
}
\end{table}
\begin{figure}[h!]
	\centering
	\begin{subfigure}[h!]{0.48\linewidth}
		\includegraphics[width=\linewidth]{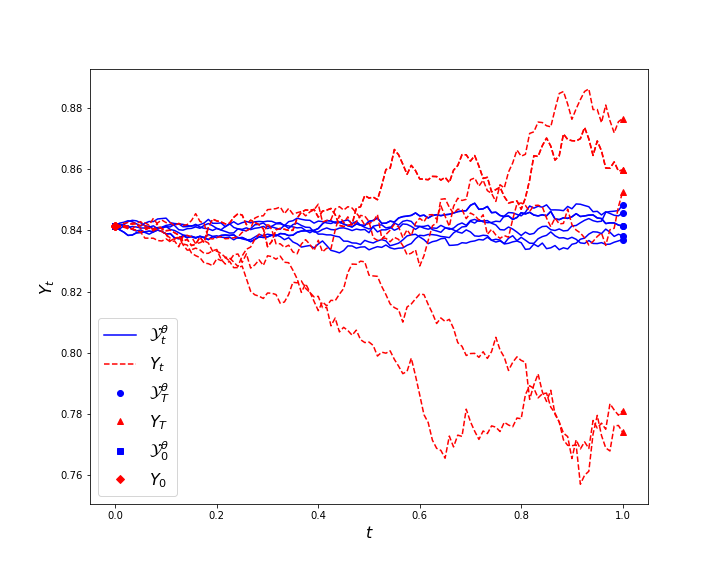}
		\caption{DBSDE $Y$ samples.}
		\label{fig6a}
	\end{subfigure}
	\begin{subfigure}[h!]{0.48\linewidth}
		\includegraphics[width=\linewidth]{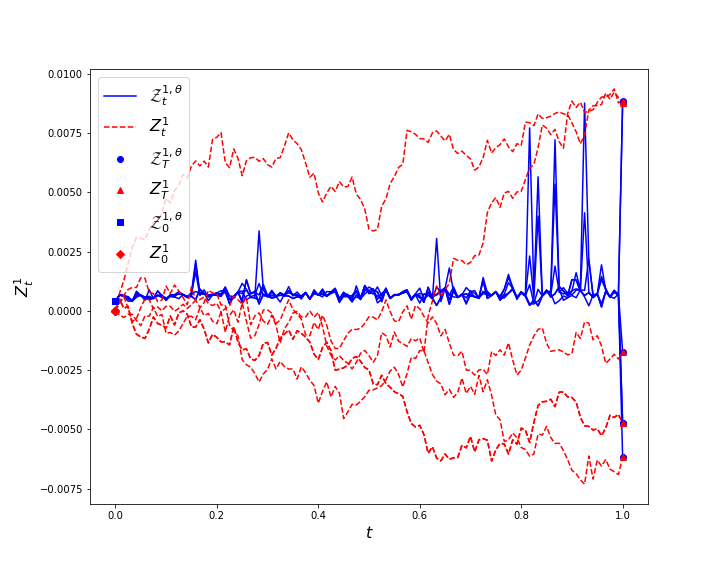}
		\caption{DBSDE $Z^1$ samples.}
		\label{fig6b}
	\end{subfigure}
	\begin{subfigure}[h!]{0.48\linewidth}
		\includegraphics[width=\linewidth]{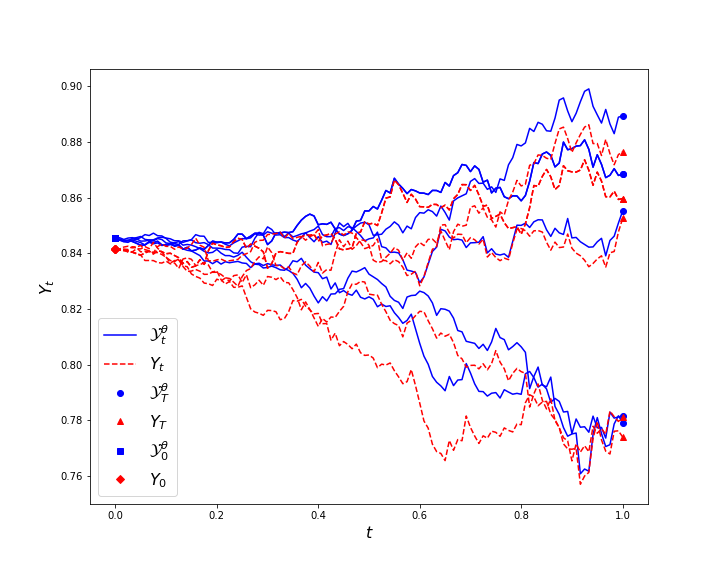}
		\caption{LDBSDE $Y$ samples.}
		\label{fig6c}
	\end{subfigure}
	\begin{subfigure}[h!]{0.48\linewidth}
		\includegraphics[width=\linewidth]{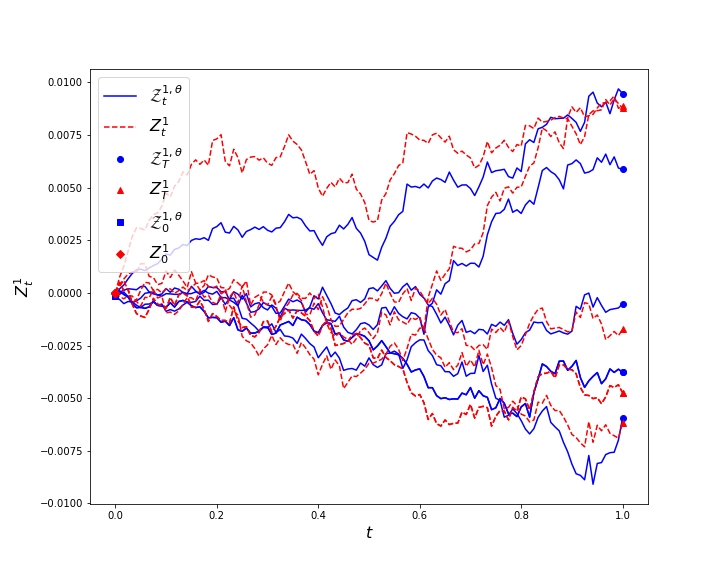}
		\caption{LDBSDE $Z^1$ samples.}
		\label{fig6d}
	\end{subfigure}	
	\begin{subfigure}[h!]{0.48\linewidth}
		\includegraphics[width=\linewidth]{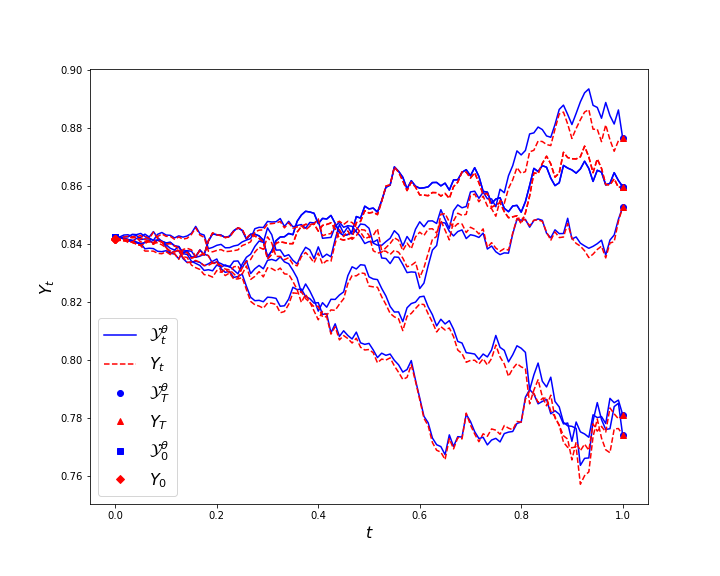}
		\caption{LaDBSDE $Y$ samples.}
		\label{fig6e}
	\end{subfigure}	
	\begin{subfigure}[h!]{0.48\linewidth}
		\includegraphics[width=\linewidth]{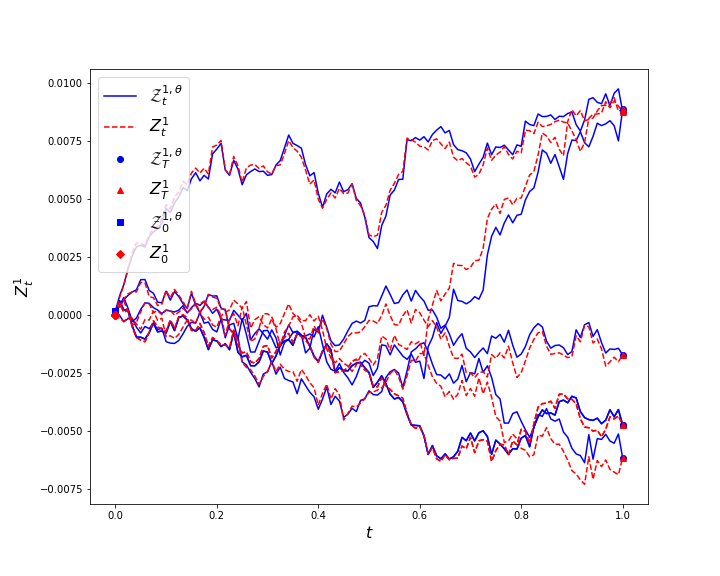}
		\caption{LaDBSDE $Z^1$ samples.}
		\label{fig6f}
	\end{subfigure}	
	 \caption{Realizations of $5$ independent paths for Example~\ref{ex2} using $d = 100$ and $N=120$. $(Y_t, Z_t^1)$ and $(\mathcal{Y}_t^{\theta},\mathcal{Z}_t^{1,\theta})$ are exact and learned solutions for $t \in [0, T]$, respectively.}
	\label{fig6}
\end{figure}
\begin{figure}[h!]
	\centering
	\begin{subfigure}[h!]{0.48\linewidth}
		\includegraphics[width=\linewidth]{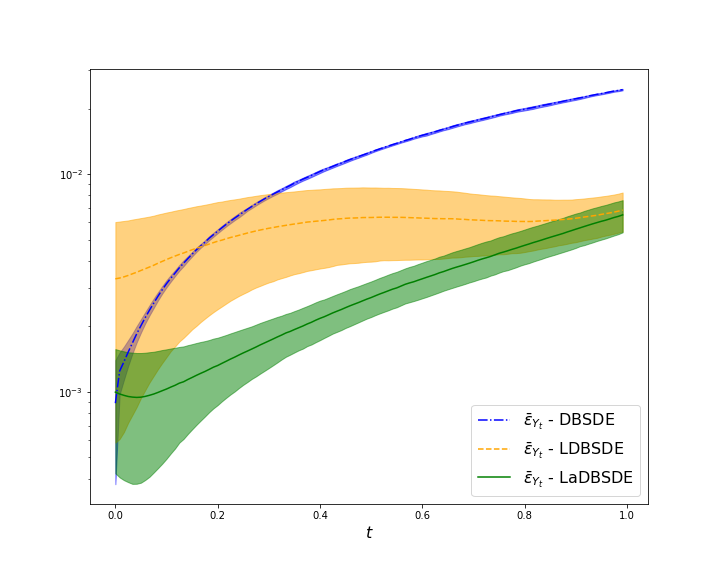}
		\caption{$Y$ process.}
		\label{fig7a}
	\end{subfigure}
	\begin{subfigure}[h!]{0.48\linewidth}
		\includegraphics[width=\linewidth]{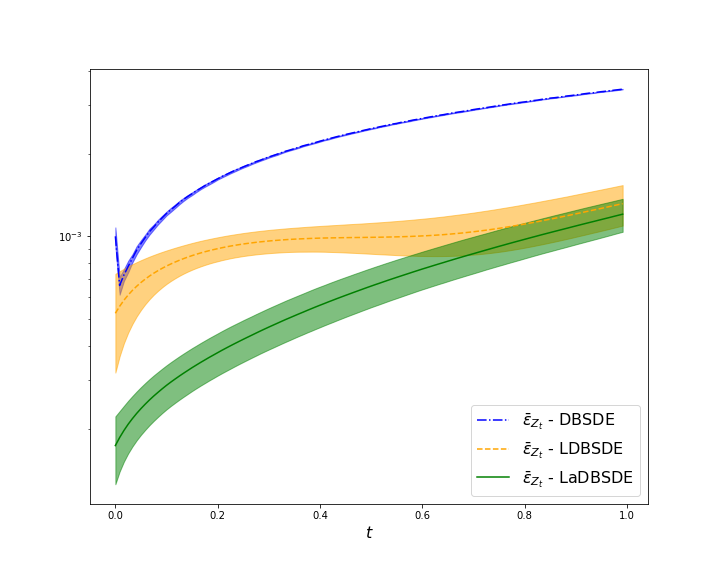}
		\caption{$Z$ process.}
		\label{fig7b}
	\end{subfigure}
	\caption{The mean regression errors  $(\bar{\epsilon}_{Y_i}, \bar{\epsilon}_{Z_i})$ at time step $t_i, i = 0, \cdots, N-1$ for Example~\ref{ex2} using $d=100$ and $N=120$. The standard deviation is given in the shaded area.}
	\label{fig7}
\end{figure}
\begin{exmp}\label{ex2}
Consider the nonlinear BSDE \cite{gobet2016linear}
\begin{equation*}
    \begin{split}
        \left\{
            \begin{array}{rcl}
                 -dY_t & = & \left( \|Z_t\|^2_{\mathbb{R}^{1\times d}} - \|\nabla \psi(t, W_t)\|^2_{\mathbb{R}^{d}} - \left( \partial_t + \frac{1}{2}\Delta \right) \psi(t, W_t) \right)\,dt -Z_t \,dW_t,\\  
       	   	    	Y_T & = & \sin\left(\| W_T\|^{2\alpha}_{\mathbb{R}^d}\right),
            \end{array}
        \right. 
    \end{split}
 \end{equation*}
\end{exmp}
where $\psi(t, W_t) = \sin\left(\left(T-t+\| W_t\|^{2}_{\mathbb{R}^d}\right)^{\alpha}\right)$. The analytic solution is 
\begin{equation*}
    \begin{split}
        \left\{
            \begin{array}{rcl}
                 Y_t & = & \sin\left(\left(T-t+\| W_t\|^{2}_{\mathbb{R}^d}\right)^{\alpha}\right),\\
   	   	    	Z_t & = & 2\alpha W_t^{\top}\cos\left(\left(T-t+\| W_t\|^{2}_{\mathbb{R}^d}\right)^{\alpha}\right)\left(T-t+\| W_t\|^{2}_{\mathbb{R}^d}\right)^{\alpha-1}.
            \end{array}
        \right.
    \end{split}
\end{equation*}

The exact solution with $d=100$, $T=1$ and $\alpha = 0.4$ is $\left(Y_0, Z_0\right) \doteq \left( 0.8415, (0, \cdots, 0)\right)$. We consider $40000$ optimization steps. We report the numerical approximation of $Y_0$ and $Z_0$ in Table~\ref{tab4} for increasing $N$. We observe comparable results for all the schemes at $t_0$.

In Figure~\ref{fig6}, we display $5$ paths of $Y$ and $Z^1$ using $N = 120$ and the regression errors in Figure~\ref{fig7}, where we see that the LaDBSDE scheme outperforms. 

For the linear and nonlinear pricing problems schemes we consider the Black-Scholes-Barenblatt type problem studied in \cite{raissi2018forward} and the problem of option pricing with different interest rates, which has been addressed in e.g., \cite{weinan2017deep, weinan2019multilevel, teng2021review, teng2022gradient}.
\begin{exmp}
Consider the Black-Scholes-Barenblatt FBSDE \cite{raissi2018forward}
\begin{equation*}
    \begin{split}
        \left\{
            \begin{array}{rcl}
                dS_t & = & \sigma S_t\, dW_t, \quad S_0 = S_0,\\ 
                -dY_t & = & -r\left(Y_t - \frac{1}{\sigma} Z_t\right)\,dt -Z_t \,dW_t,\\  
   	   	    	Y_T & = & \|S_T\|^2_{\mathbb{R}^d},
            \end{array}
        \right.
    \end{split}
\end{equation*}
\label{ex3} 
\end{exmp}
 The analytic solution is 
\begin{equation*}
    \begin{split}
        \left\{
            \begin{array}{rcl}
                 Y_t & = & \exp\left(\left(r + \sigma^2 \right) (T-t)\right)\|S_t\|^2_{\mathbb{R}^d},\\
   	   	    	Z_t & = & 2\sigma \exp\left(\left(r + \sigma^2 \right) (T-t)\right)S_t^2.
            \end{array}
        \right.
    \end{split}
\end{equation*}
\begin{figure}[h!]
	\centering
	\begin{subfigure}[h!]{0.48\linewidth}
		\includegraphics[width=\linewidth]{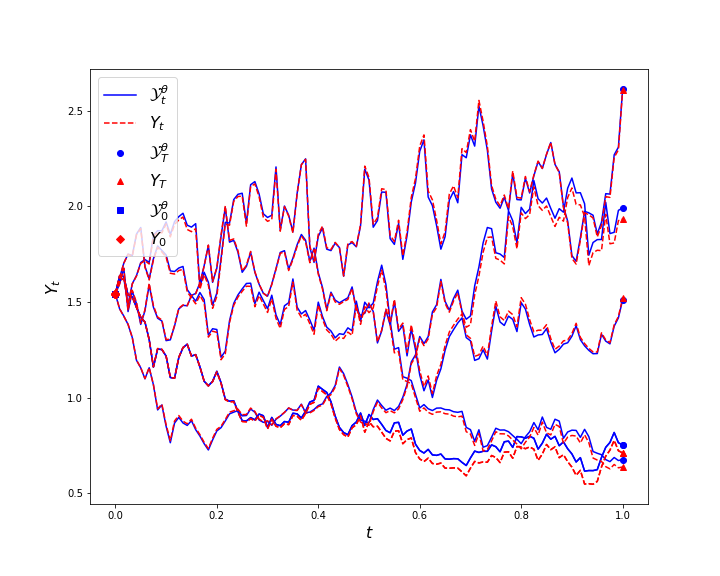}
		\caption{DBSDE $Y$ samples.}
		\label{fig8a}
	\end{subfigure}
	\begin{subfigure}[h!]{0.48\linewidth}
		\includegraphics[width=\linewidth]{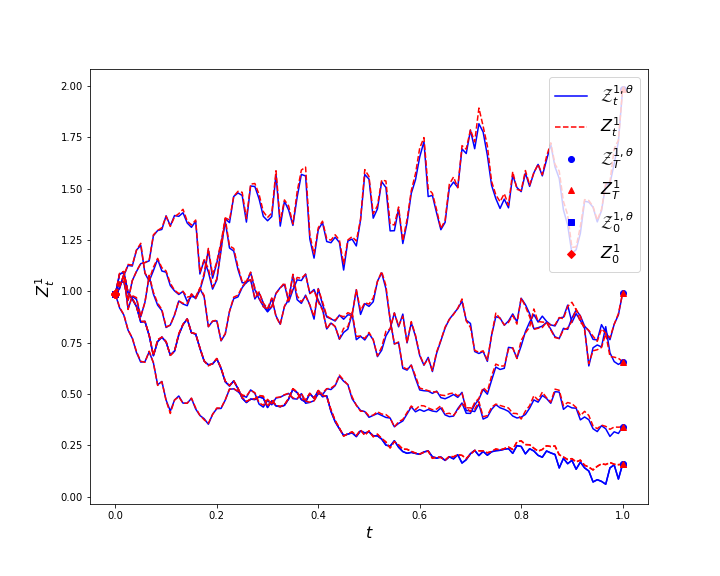}
		\caption{DBSDE $Z^1$ samples.}
		\label{fig8b}
	\end{subfigure}
	\begin{subfigure}[h!]{0.48\linewidth}
		\includegraphics[width=\linewidth]{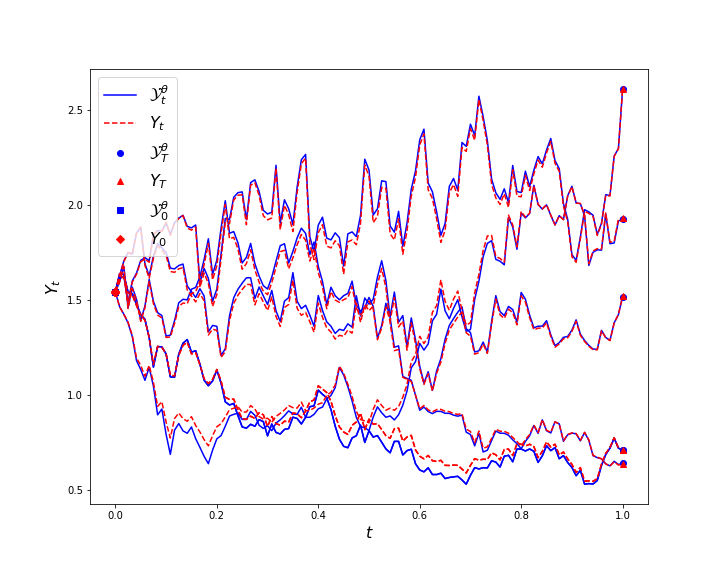}
		\caption{LDBSDE $Y$ samples.}
		\label{fig8c}
	\end{subfigure}
	\begin{subfigure}[h!]{0.48\linewidth}
		\includegraphics[width=\linewidth]{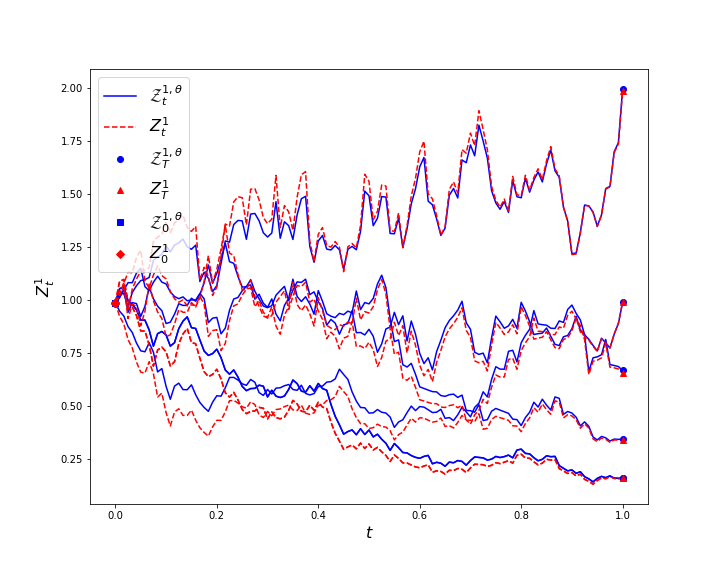}
		\caption{LDBSDE $Z^1$ samples.}
		\label{fig8d}
	\end{subfigure}
	\begin{subfigure}[h!]{0.48\linewidth}
		\includegraphics[width=\linewidth]{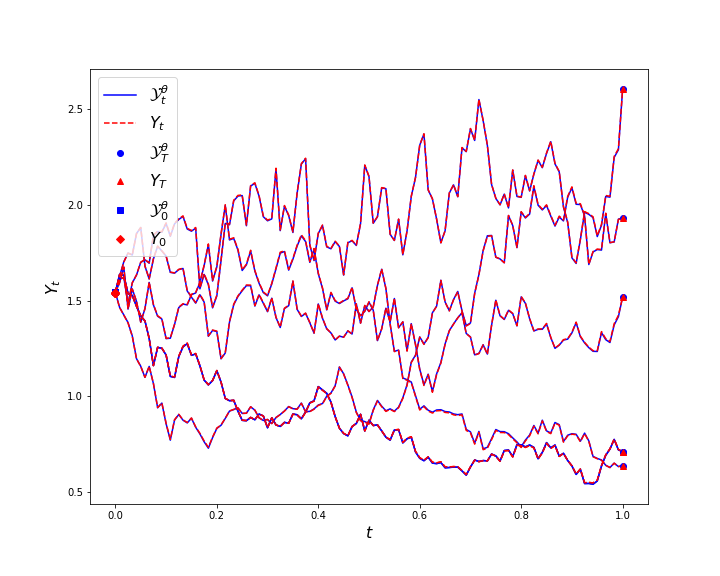}
		\caption{LaDBSDE $Y$ samples.}
		\label{fig8e}
	\end{subfigure}
	\begin{subfigure}[h!]{0.48\linewidth}
		\includegraphics[width=\linewidth]{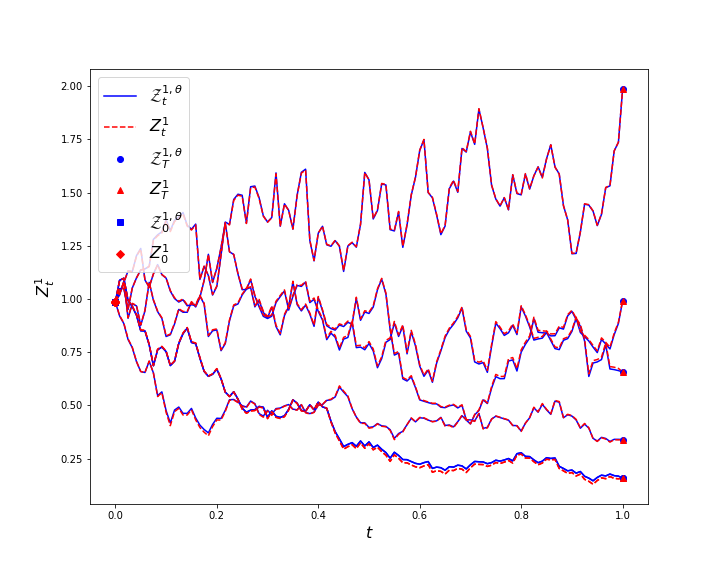}
		\caption{LaDBSDE $Z^1$ samples.}
		\label{fig8f}
	\end{subfigure}
	 \caption{Realizations of $5$ independent paths for Example~\ref{ex3} using $d = 2$ and $N=120$. $(Y_t, Z_t^1)$ and $(\mathcal{Y}_t^{\theta},\mathcal{Z}_t^{1,\theta})$ are exact and learned solutions for $t \in [0, T]$, respectively.}
	\label{fig8}
\end{figure}
\begin{figure}[h!]
	\centering
	\begin{subfigure}[h!]{0.48\linewidth}
		\includegraphics[width=\linewidth]{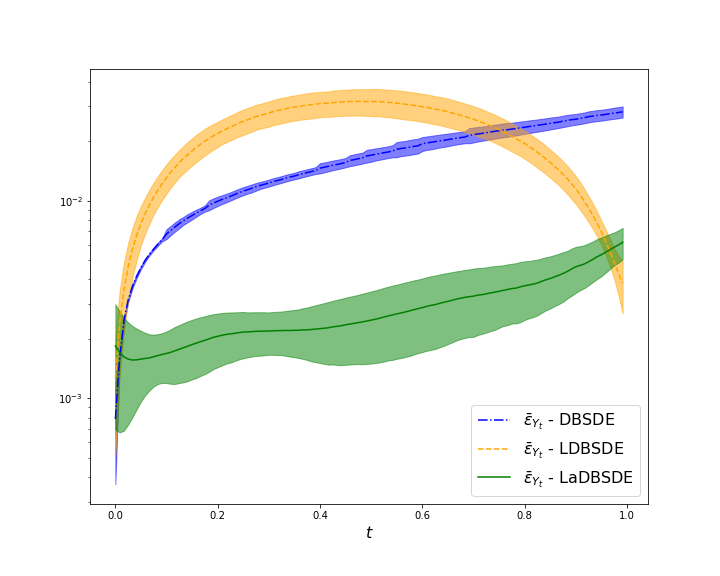}
		\caption{$Y$ process.}
		\label{fig9a}
	\end{subfigure}
	\begin{subfigure}[h!]{0.48\linewidth}
		\includegraphics[width=\linewidth]{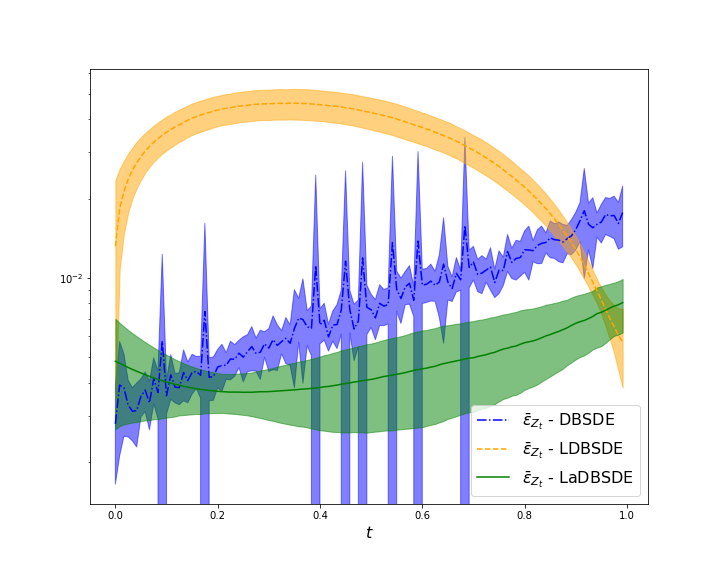}
		\caption{$Z$ process.}
		\label{fig9b}
	\end{subfigure}
	\caption{The mean regression errors  $(\bar{\epsilon}_{Y_i}, \bar{\epsilon}_{Z_i})$ at time step $t_i, i = 0, \cdots, N-1$ for Example~\ref{ex3} using $d=2$ and $N=120$. The standard deviation is given in the shaded area.}
	\label{fig9}
\end{figure}
We use $T=1$, $r = 0.05$, $\sigma = 0.4$ and $S_0 = (1, 0.5, \cdots, 1, 0.5) \in \mathbb{R}^d$. We start in the case of $d=2$. The exact solution is $\left(Y_0, Z_0\right) \doteq \left( 1.5421, (0.9869, 0.2467)\right)$. Using $40000$ optimization steps and $N = 120$, the numerical approximation of $Y_0$ and $Z_0$ is given in Table~\ref{tab5}. 
\begin{table}[h!]
{\footnotesize
\begin{center}
  \begin{tabular}{| c | c | c |}
  \hline
   Scheme &  $\bar{\epsilon}_{Y_0}$ (Std. Dev.) & $\bar{\epsilon}_{Z_0}$ (Std. Dev.)\\ \hline
   DBSDE & 7.87e-4 (4.19e-4) & 2.80e-3 (1.15e-3) \\ \hline
   LDBSDE & 1.26e-3 (7.50e-4) & 1.32e-2 (1.03e-2) \\ \hline
   LaDBSDE & 1.84e-3 (1.15e-3) & 4.83e-3 (2.17e-3) \\ \hline
   \end{tabular}
  \end{center}
\caption{The mean absolute errors of $Y_0$ and $Z_0$ for Example~\ref{ex3} using $d = 2$ and $N = 120$. The standard deviation is given in parenthesis.}
\label{tab5}  
}
\end{table}
The DBSDE scheme gives smaller errors at $t_0$ compared to the schemes LDBSDE and LaDBSDE. 

However, our scheme gives the best approximations for $t>t_0$. This can be observed in Figures~\ref{fig8} and~\ref{fig9}, where $5$ paths of $Y$ and $Z^1$ and the regression errors are displayed, respectively.

\begin{figure}[h!]
	\centering
	\begin{subfigure}[h!]{0.48\linewidth}
		\includegraphics[width=\linewidth]{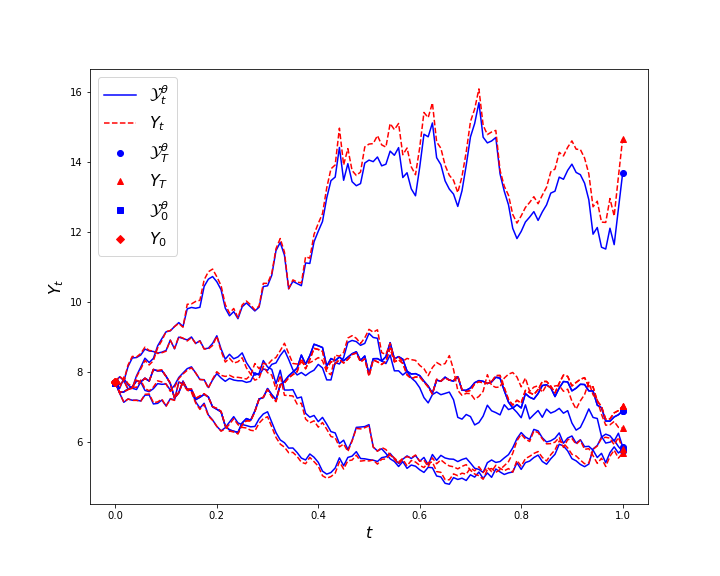}
		\caption{DBSDE $Y$ samples.}
		\label{fig10a}
	\end{subfigure}
	\begin{subfigure}[h!]{0.48\linewidth}
		\includegraphics[width=\linewidth]{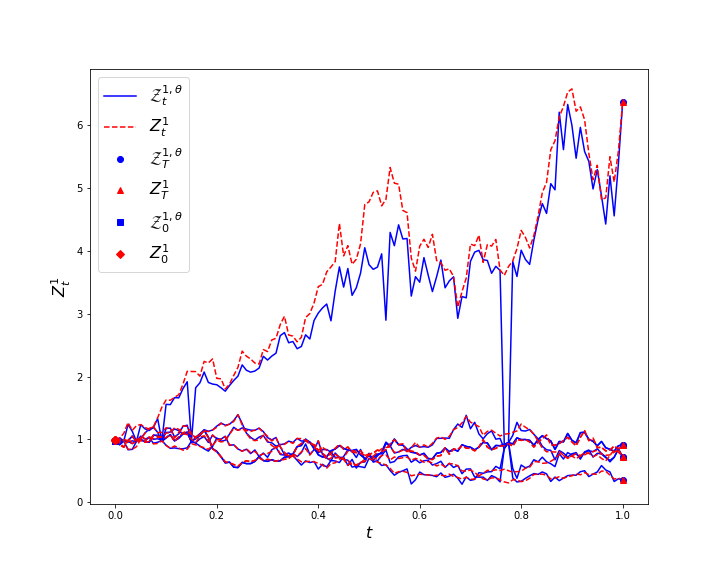}
		\caption{DBSDE $Z^1$ samples.}
		\label{fig10b}
	\end{subfigure}
	\begin{subfigure}[h!]{0.48\linewidth}
		\includegraphics[width=\linewidth]{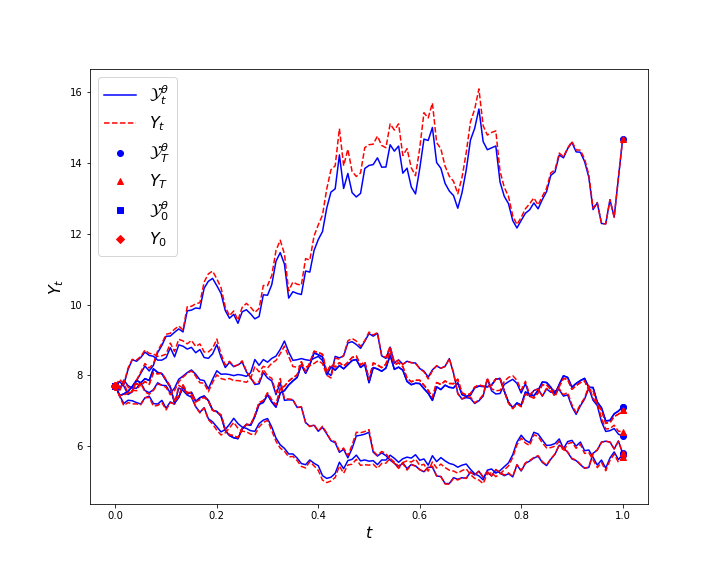}
		\caption{LDBSDE $Y$ samples.}
		\label{fig10c}
	\end{subfigure}
	\begin{subfigure}[h!]{0.48\linewidth}
		\includegraphics[width=\linewidth]{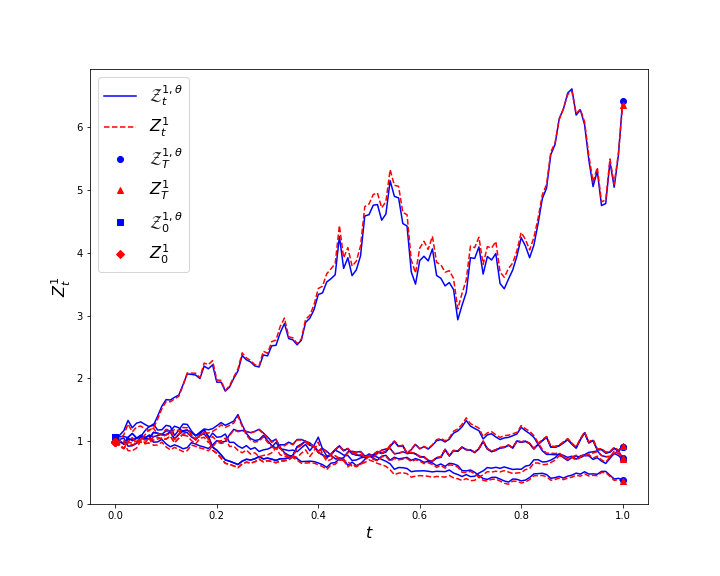}
		\caption{LDBSDE $Z^1$ samples.}
		\label{fig10d}
	\end{subfigure}
	\begin{subfigure}[h!]{0.48\linewidth}
		\includegraphics[width=\linewidth]{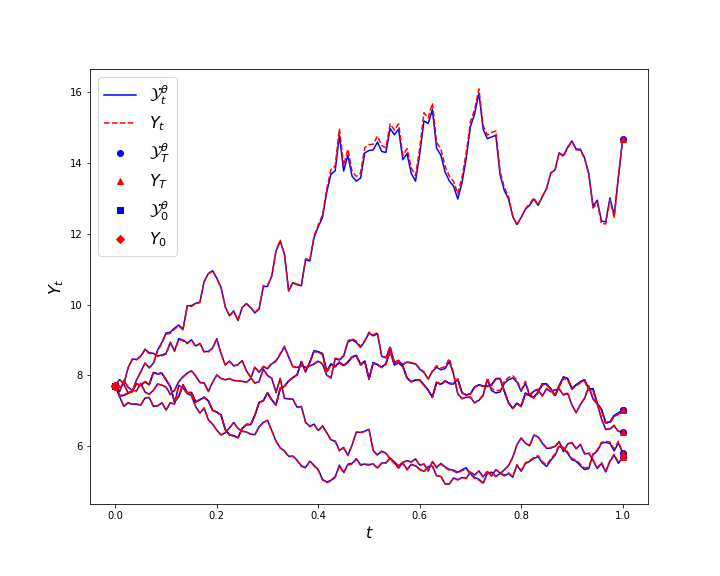}
		\caption{LaDBSDE $Y$ samples.}
		\label{fig10e}
	\end{subfigure}
	\begin{subfigure}[h!]{0.48\linewidth}
		\includegraphics[width=\linewidth]{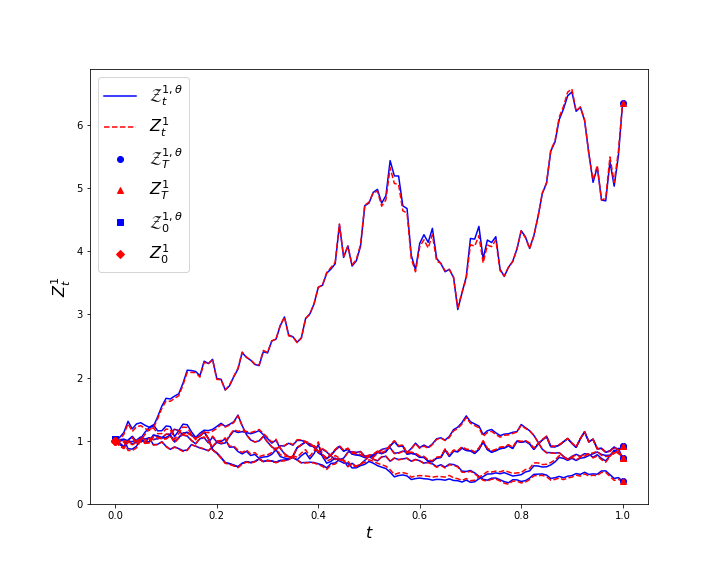}
		\caption{LaDBSDE $Z^1$ samples.}
		\label{fig10f}
	\end{subfigure}
	 \caption{Realizations of $5$ independent paths for Example~\ref{ex3} using $d = 10$ and $N=120$. $(Y_t, Z_t^1)$ and $(\mathcal{Y}_t^{\theta},\mathcal{Z}_t^{1,\theta})$ are exact and learned solutions for $t \in [0, T]$, respectively.}
	\label{fig10}
\end{figure}

\begin{figure}[h!]
	\centering
	\begin{subfigure}[h!]{0.48\linewidth}
		\includegraphics[width=\linewidth]{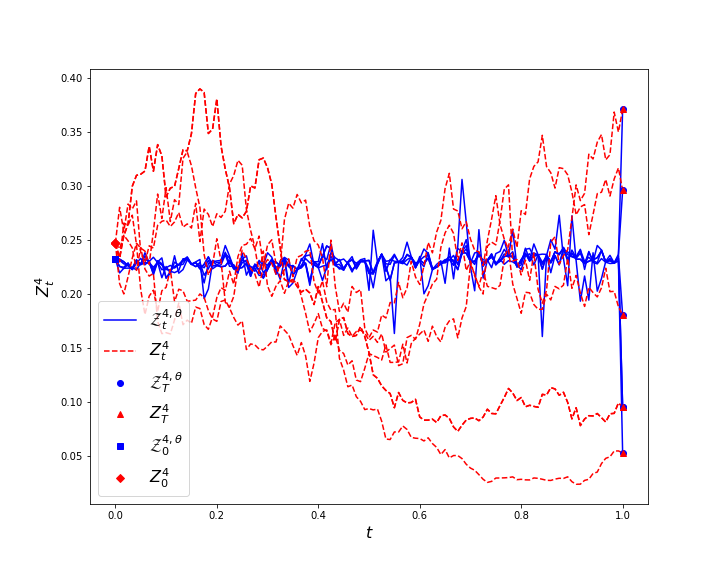}
		\caption{DBSDE $Z^4$ samples.}
		\label{fig11a}
	\end{subfigure}
	\begin{subfigure}[h!]{0.48\linewidth}
		\includegraphics[width=\linewidth]{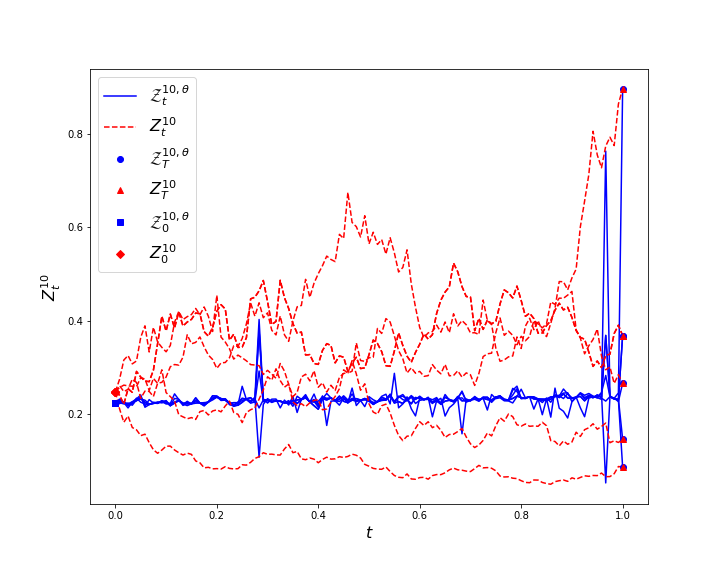}
		\caption{DBSDE $Z^{10}$ samples.}
		\label{fig11b}
	\end{subfigure}
	\begin{subfigure}[h!]{0.48\linewidth}
		\includegraphics[width=\linewidth]{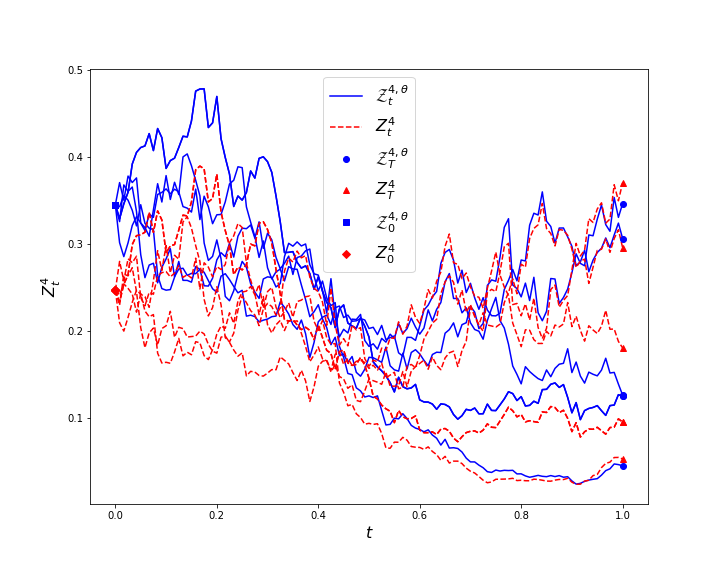}
		\caption{LDBSDE $Z^4$ samples.}
		\label{fig11c}
	\end{subfigure}
	\begin{subfigure}[h!]{0.48\linewidth}
		\includegraphics[width=\linewidth]{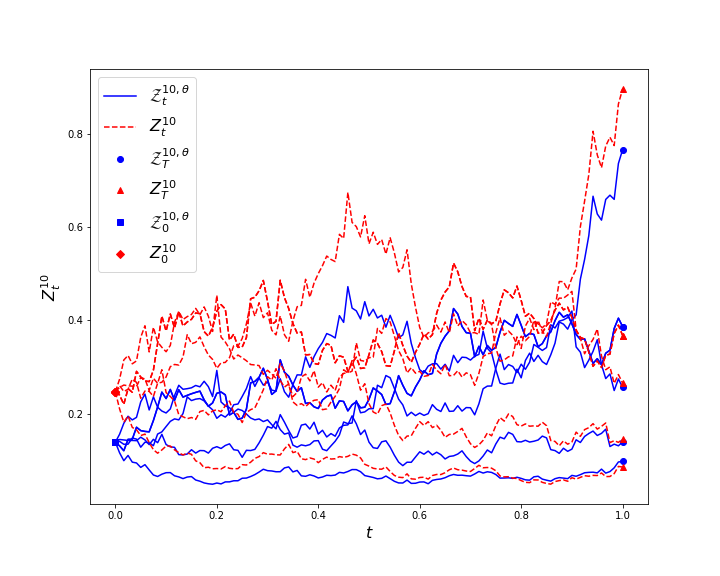}
		\caption{LDBSDE $Z^{10}$ samples.}
		\label{fig11d}
	\end{subfigure}
	\begin{subfigure}[h!]{0.48\linewidth}
		\includegraphics[width=\linewidth]{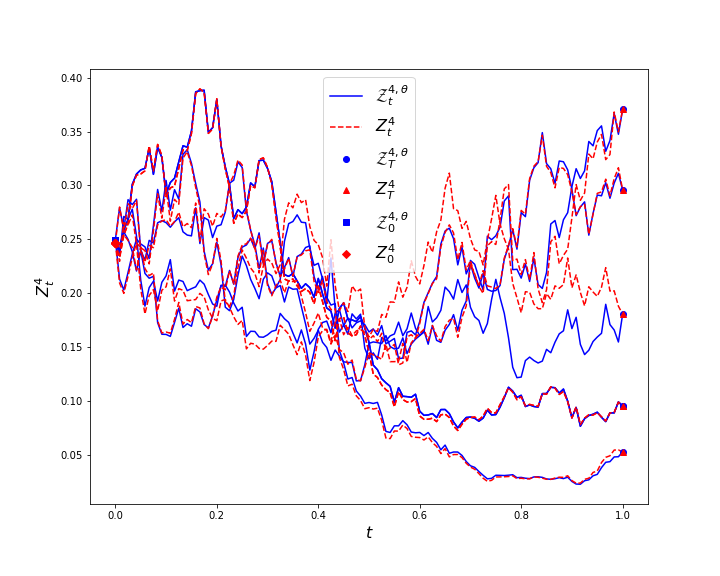}
		\caption{LaDBSDE $Z^4$ samples.}
		\label{fig11e}
	\end{subfigure}
	\begin{subfigure}[h!]{0.48\linewidth}
		\includegraphics[width=\linewidth]{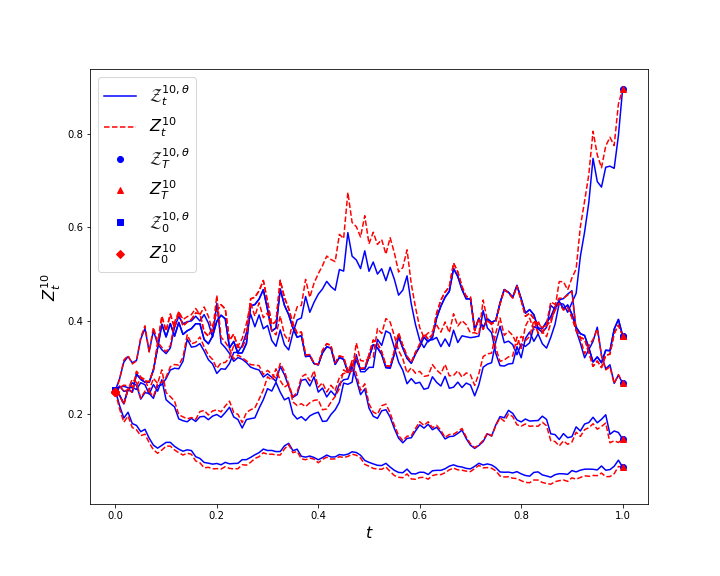}
		\caption{LaDBSDE $Z^{10}$ samples.}
		\label{fig11f}
	\end{subfigure}
	 \caption{Realizations of $5$ independent paths for Example~\ref{ex3} using $d = 10$ and $N=120$. $(Z_t^4, Z_t^{10})$ and $(\mathcal{Z}_t^{4, \theta},\mathcal{Z}_t^{10,\theta})$ are exact and learned solutions for $t \in [0, T]$, respectively.}
	\label{fig11}
\end{figure}

\begin{figure}[h!]
	\centering
	\begin{subfigure}[h!]{0.48\linewidth}
		\includegraphics[width=\linewidth]{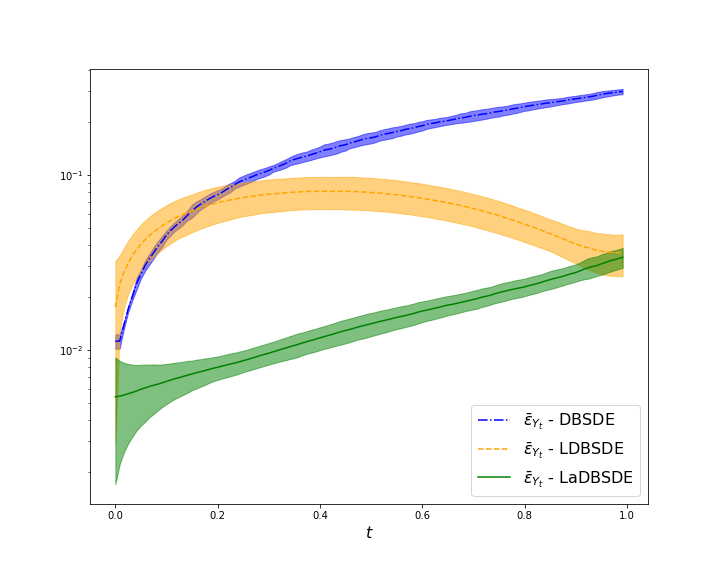}
		\caption{$Y$ process.}
		\label{fig12a}
	\end{subfigure}
	\begin{subfigure}[h!]{0.48\linewidth}
		\includegraphics[width=\linewidth]{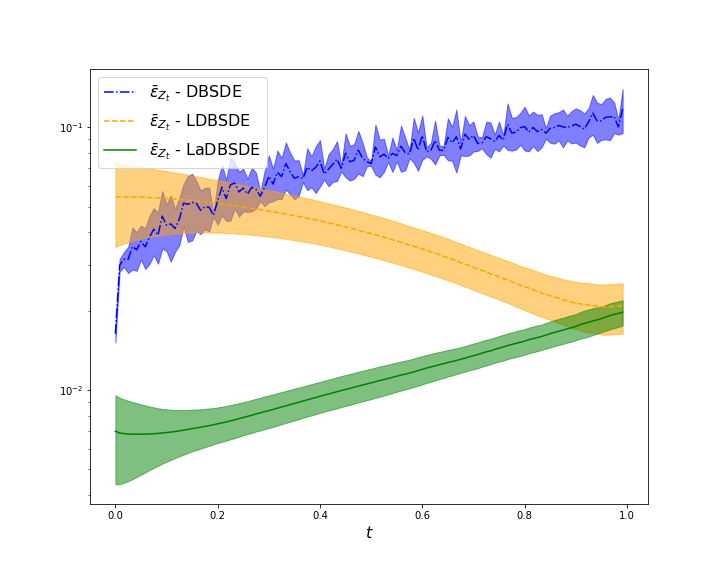}
		\caption{$Z$ process.}
		\label{fig12b}
	\end{subfigure}
	\caption{The mean regression errors  $(\bar{\epsilon}_{Y_i}, \bar{\epsilon}_{Z_i})$ at time step $t_i, i = 0, \cdots, N-1$ for Example~\ref{ex3} using $d=10$ and $N=120$. The standard deviation is given in the shaded area.}
	\label{fig12}
\end{figure}

Now we increase the dimension by setting $d=10$. The exact solution is $(Y_0, Z_0) \doteq \left( 7.7105, (0.9869, 0.2467, \cdots, 0.9869, 0.2467)\right).$ The numerical approximations of $Y_0$ and $Z_0$ using $40000$ optimization steps and $N = 120$ are reported in Table~\ref{tab6}.
\begin{table}[h!]
{\footnotesize
\begin{center}
  \begin{tabular}{| c | c | c |}
  \hline
   Scheme &  $\bar{\epsilon}_{Y_0}$ (Std. Dev.) & $\bar{\epsilon}_{Z_0}$ (Std. Dev.)\\ \hline
   DBSDE & 1.12e-2 (1.09e-3) & 1.64e-2 (1.22e-3) \\ \hline
   LDBSDE & 1.76e-2 (1.46e-2) & 5.44e-2 (1.94e-2) \\ \hline
   LaDBSDE & 5.39e-3 (3.68e-3) & 6.98e-3 (2.59e-3) \\ \hline
   \end{tabular}
  \end{center}
\caption{The mean absolute errors of $Y_0$ and $Z_0$ for Example~\ref{ex3} using $d = 10$ and $N = 120$. The standard deviation is given in parenthesis.}
\label{tab6}  
}
\end{table}
Our scheme gives the smallest errors. Using $5$ paths of $Y$ and $Z^1$, we compare the approximations for the entire time domain in Figure~\ref{fig10}. Note that the approximation quality of each component in $Z$ may be different. To show this we display the approximations of $Z^4$ and $Z^{10}$ in Figure~\ref{fig11}. The DBSDE scheme fails to perform well for each component of process $Z$, whereas our scheme maintains its robustness. Furthermore, the LaDBSDE scheme provides the smallest regression errors as shown in Figure~\ref{fig12}.

\begin{figure}[htb!]
	\centering
	\begin{subfigure}[h!]{0.48\linewidth}
		\includegraphics[width=\linewidth]{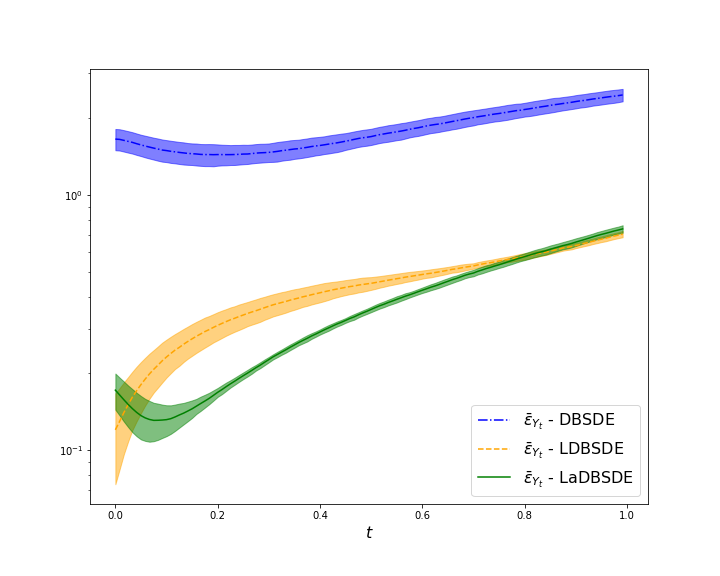}
		\caption{$Y$ process.}
		\label{fig13a}
	\end{subfigure}
	\begin{subfigure}[h!]{0.48\linewidth}
		\includegraphics[width=\linewidth]{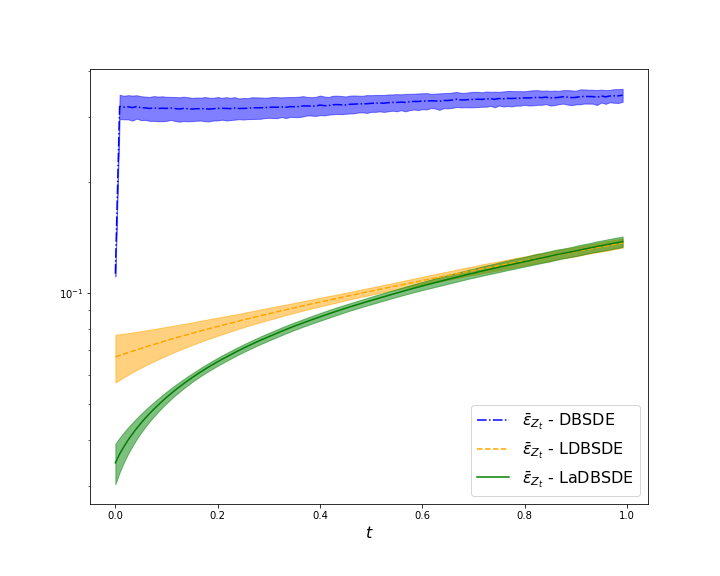}
		\caption{$Z$ process.}
		\label{fig13b}
	\end{subfigure}
	\caption{The mean regression errors  $(\bar{\epsilon}_{Y_i}, \bar{\epsilon}_{Z_i})$ at time step $t_i, i = 0, \cdots, N-1$ for Example~\ref{ex3} using $d=50$ and $N=120$. The standard deviation is given in the shaded area.}
	\label{fig13}
\end{figure}

\begin{figure}[htb!]
	\centering
	\begin{subfigure}[h!]{0.32\linewidth}
		\includegraphics[width=\linewidth]{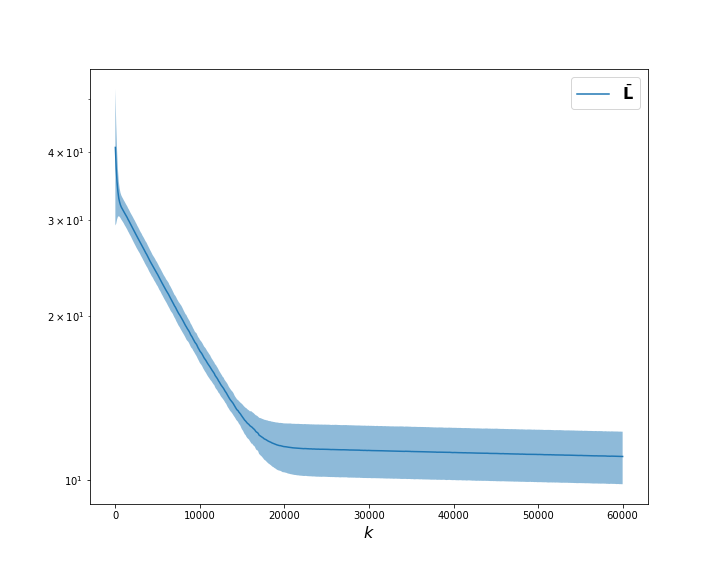}
		\caption{DBSDE validation loss.}
		\label{fig14a}
	\end{subfigure}
	\begin{subfigure}[h!]{0.32\linewidth}
		\includegraphics[width=\linewidth]{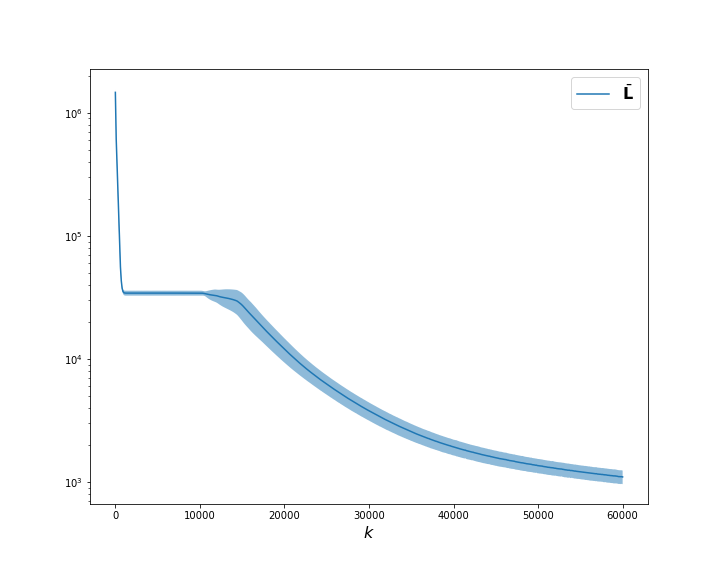}
		\caption{LDBSDE validation loss.}
		\label{fig14b}
	\end{subfigure}
	\begin{subfigure}[h!]{0.32\linewidth}
		\includegraphics[width=\linewidth]{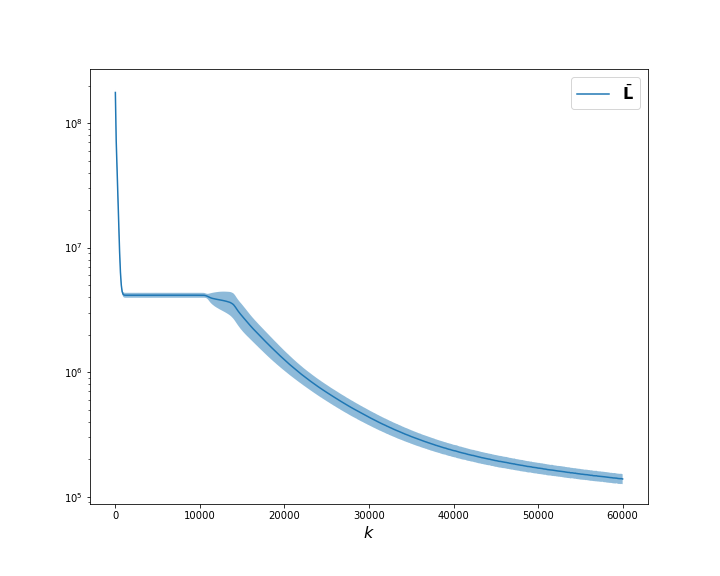}
		\caption{LaDBSDE validation loss.}
		\label{fig14c}
	\end{subfigure}	
	\caption{The mean loss values  $\bar{\mathbf{L}}$ for Example~\ref{ex3} using $d=50$ and $N=120$. The standard deviation is given in the shaded area.}
	\label{fig14}
\end{figure}

\begin{figure}[t!]
	\centering
	\begin{subfigure}[h!]{0.48\linewidth}
		\includegraphics[width=\linewidth]{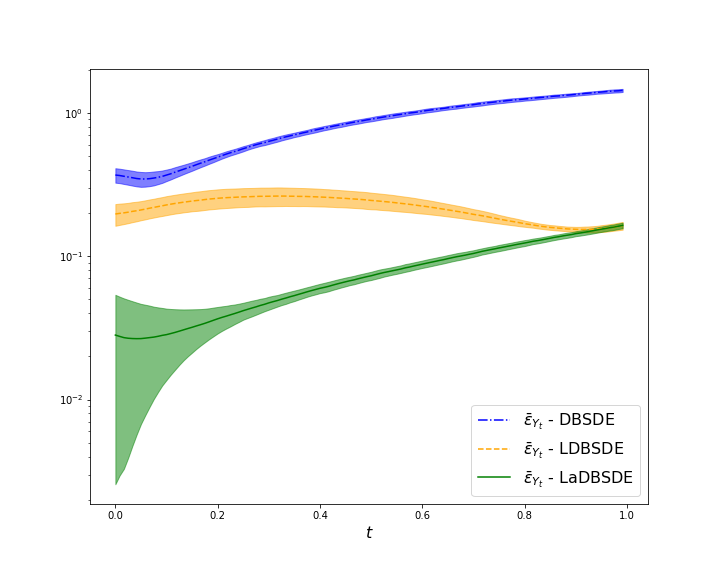}
		\caption{$Y$ process.}
		\label{fig15a}
	\end{subfigure}
	\begin{subfigure}[h!]{0.48\linewidth}
		\includegraphics[width=\linewidth]{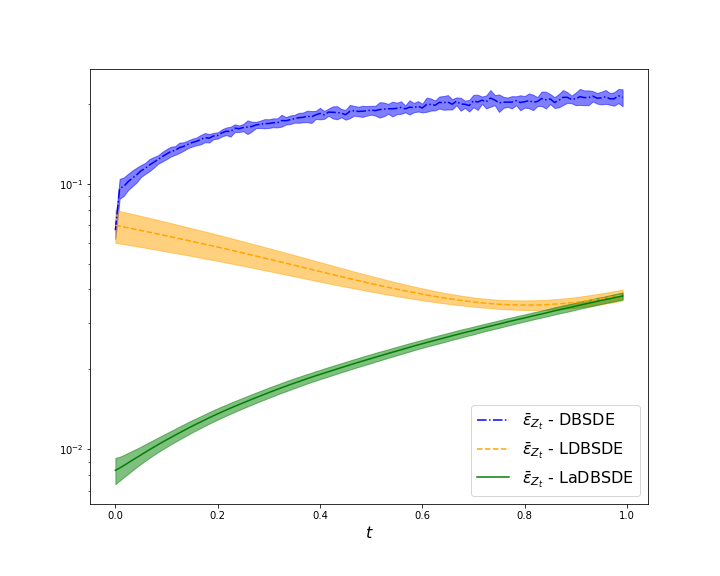}
		\caption{$Z$ process.}
		\label{fig15b}
	\end{subfigure}
	\caption{The mean regression errors  $(\bar{\epsilon}_{Y_i}, \bar{\epsilon}_{Z_i})$ at time step $t_i, i = 0, \cdots, N-1$ for Example~\ref{ex3} using $d=50$ and $N=120$. The standard deviation is given in the shaded area.}
	\label{fig15}
\end{figure}

We further increase the dimension by setting $d=50$. The exact solution  is $(Y_0, Z_0) \doteq \left( 38.5524, (0.9869, 0.2467, \cdots, 0.9869, 0.2467)\right).$ We use $60000$ optimization steps. For $N=120$, the numerical approximations of $Y_0$ and $Z_0$ are given in Table~\ref{tab7}, we see that the schemes LDBSDE and LaDBSDE perform similarly, and better than the DBSDE scheme.
\begin{table}[h!]
{\footnotesize
\begin{center}
  \begin{tabular}{| c | c | c |}
  \hline
   Scheme &  $\bar{\epsilon}_{Y_0}$ (Std. Dev.) & $\bar{\epsilon}_{Z_0}$ (Std. Dev.)\\ \hline
   DBSDE & 1.66e+0 (1.59e-1) & 1.13e-1 (2.00e-3) \\ \hline
   LDBSDE & 1.20e-1 (4.67e-2) & 6.71e-2 (9.97e-3) \\ \hline
   LaDBSDE & 1.72e-1 (2.78e-2) & 3.46e-2 (4.33e-3) \\ \hline
   \end{tabular}
  \end{center}
\caption{The mean absolute errors of $Y_0$ and $Z_0$ for Example~\ref{ex3} using $d = 50$ and $N = 120$. The standard deviation is given in parenthesis.}
\label{tab7}  
}
\end{table}
Furthermore, the smallest regression errors are provided by the LaDBSDE scheme as displayed in Figure~\ref{fig13}.
Note that the results can be further improved as it can be seen from the validation plots of the mean loss value $\bar{\mathbf{L}} = \frac{1}{10}\sum_{i=1}^{10} \mathbf{L}_i$ in Figure~\ref{fig14}.
To do that, for the first $30000$ optimization steps we use the learning rate $\gamma_0$ and apply the learning rate decay approach for the next $30000$ optimization steps. The numerical approximations of $Y_0$ and $Z_0$ for $N = 120$ are given in Table~\ref{tab8} and the regression errors in Figure~\ref{fig15}. We see that the LaDBSDE scheme outperforms always.
\begin{table}[t!]
{\footnotesize
\begin{center}
  \begin{tabular}{| c | c | c |}
  \hline
   Scheme &  $\bar{\epsilon}_{Y_0}$ (Std. Dev.) & $\bar{\epsilon}_{Z_0}$ (Std. Dev.)\\ \hline
   DBSDE & 3.68e-1 (4.26e-2) & 6.69e-2 (4.90e-3) \\ \hline
   LDBSDE & 1.97e-1 (3.40e-2) & 6.97e-2 (9.85e-3) \\ \hline
   LaDBSDE & 2.82e-2 (2.56e-2) & 8.35e-3 (9.51e-4) \\ \hline
   \end{tabular}
  \end{center}
\caption{The mean absolute errors of $Y_0$ and $Z_0$ for Example~\ref{ex3} using $d = 50$ and $N = 120$. The standard deviation is given in parenthesis.}
\label{tab8}  
}
\end{table}

\begin{figure}[htb!]
	\centering
	\begin{subfigure}[h!]{0.48\linewidth}
		\includegraphics[width=\linewidth]{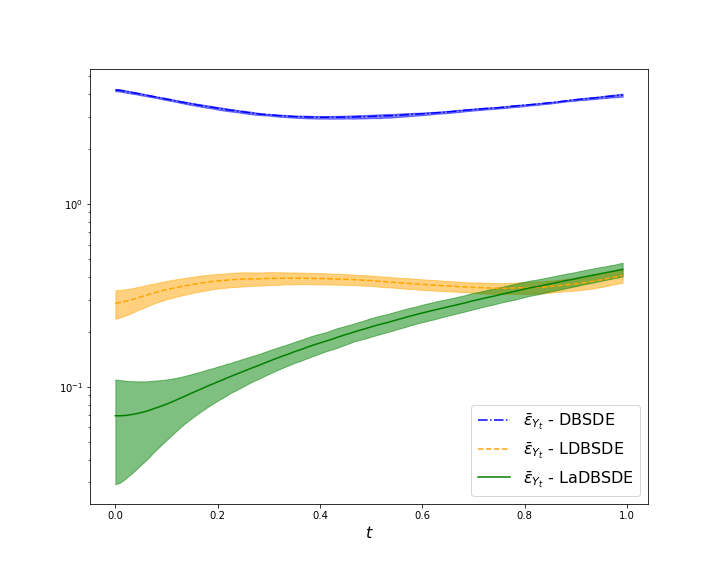}
		\caption{$Y$ process.}
		\label{fig16a}
	\end{subfigure}
	\begin{subfigure}[h!]{0.48\linewidth}
		\includegraphics[width=\linewidth]{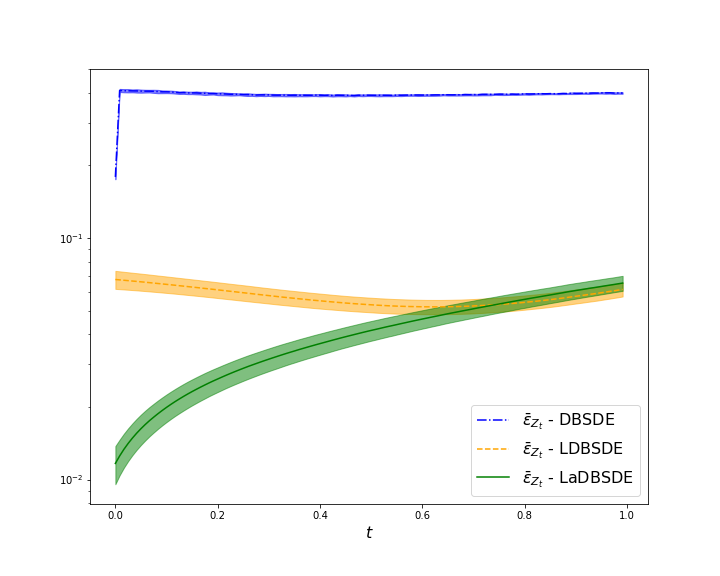}
		\caption{$Z$ process.}
		\label{fig16b}
	\end{subfigure}
	\caption{The mean regression errors  $(\bar{\epsilon}_{Y_i}, \bar{\epsilon}_{Z_i})$ at time step $t_i, i = 0, \cdots, N-1$ for Example~\ref{ex3} using $d=100$ and $N=120$. The standard deviation is given in the shaded area.}
	\label{fig16}
\end{figure}

Finally, we consider $d=100$ with $(Y_0, Z_0) \doteq \left( 77.1049, (0.9869, 0.2467, \cdots, 0.9869, 0.2467)\right).$ We use the same technique for the learning approach as that in $d=50$ in order to improve the results. The numerical approximation of $Y_0$ and $Z_0$ for an increasing $N$ is reported in Table~\ref{tab9}.
\begin{table}[h!]
{\footnotesize
\begin{center}
  \begin{tabular}{| c | c | c | c | c |}
  \hline
   \multirow{3}{*}{Scheme} & N = 60 & N = 80 & N = 100 & N = 120\\ \
    &  $\bar{\epsilon}_{Y_0}$ (Std. Dev.) & $\bar{\epsilon}_{Y_0}$ (Std. Dev.) & $\bar{\epsilon}_{Y_0}$ (Std. Dev.) & $\bar{\epsilon}_{Y_0}$ (Std. Dev.)\\ 
    &  $\bar{\epsilon}_{Z_0}$ (Std. Dev.) & $\bar{\epsilon}_{Z_0}$ (Std. Dev.) & $\bar{\epsilon}_{Z_0}$ (Std. Dev.) & $\bar{\epsilon}_{Z_0}$ (Std. Dev.)\\ \hline 
    \multirow{2}{*}{DBSDE} & 3.25e+0 (7.40e-2) & 3.65e+0 (7.23e-2) & 3.96e+0 (6.82e-2) & 4.20e+0 (6.32e-2) \\
    & 1.79e-1 (2.75e-3) & 1.82e-1 (2.99e-3) & 1.83e-1 (3.19e-3) & 1.79e-1 (4.08e-3) \\ \hline
   \multirow{2}{*}{LDBSDE} & 3.10e-1 (4.24e-2) & 2.73e-1 (4.75e-2) & 3.07e-1 (7.04e-2) & 2.86e-1 (5.11e-2) \\
    & 5.77e-2 (4.52e-3) & 6.13e-2 (5.39e-3) & 6.70e-2 (6.94e-3) & 6.74e-2 (5.77e-3) \\ \hline
   \multirow{2}{*}{LaDBSDE} & 7.14e-2 (3.97e-2) & 5.66e-2 (4.69e-2) & 3.88e-2 (3.99e-2) & 6.95e-2 (4.02e-2)\\
   & 1.12e-2 (2.58e-3) & 1.21e-2 (1.49e-3) & 1.13e-2 (3.23e-3) & 1.17e-2 (2.10e-3) \\ \hline
   \end{tabular}
  \end{center}
 \caption{The mean absolute errors of $Y_0$ and $Z_0$ for Example~\ref{ex3} using $d = 100$. The standard deviation is given in parenthesis.}
\label{tab9} 
}
\end{table}
The same conclusion can be drawn that the LaDBSDE scheme outperforms. More precisely, the relative error for the DBSDE method with $N = 120$ for $Y_0$ is $5.45\%$ and $17.90\%$ for $Z_0$. The LDBSDE scheme achieves $0.37\%$ and $6.74\%$ respectively, while the LaDBSDE method gives $0.09\%$ and $1.17\%$. The regression errors are displayed in Figure~\ref{fig16}.

\begin{exmp}
Consider the nonlinear pricing with different interest rates \cite{bergman1995option}
\begin{equation*}
    \begin{split}
        \left\{
            \begin{array}{rcl}
                dS_t & = & \mu S_t\,dt + \sigma S_t\, dW_t,  \quad S_0 = S_0,\\ 
                -dY_t & = & \left(-R^lY_t - \frac{\mu - R^l}{\sigma} \sum_{i=1}^{d}Z_t^i+ \left( R^b - R^l\right) \max \left( \frac{1}{\sigma} \sum_{i=1}^{d}Z_t^i - Y_t , 0 \right)\right)\,dt -Z_t \,dW_t,\\  
   	   	    	Y_T & = & \max \left(\max_{d = 1, \cdots, D} (S_{T,d} - K_1, 0\right)-2\max \left(\max_{d = 1, \cdots, D} (S_{T,d} - K_2, 0\right),
            \end{array}
        \right.
    \end{split}
\end{equation*}
\label{ex4}
\end{exmp}
where $S_t = (S_t^1, S_t^2, \cdots, S_t^d)^{\top}$. The benchmark value with $T=0.5$, $\mu = 0.06$, $\sigma = 0.2$, $R^l = 0.04$, $R^b = 0.06$, $K_1 = 120$, $K_2 = 150$ and $S_0 = 100$ is $Y_0 \doteq 21.2988$, which is computed using the multilevel Monte Carlo with 7 Picard iterations approach \cite{weinan2019multilevel}. We use $30000$ optimization steps, and show numerical approximation for $Y_0$ (the reference results for $Z_0$ are not available) for an increasing $N$ in Table~\ref{tab10}. The approximations by all the schemes are comparable.
\begin{table}[h!]
{\footnotesize
\begin{center}
  \begin{tabular}{| c | c | c | c | c |}
  \hline
     \multirow{2}{*}{Scheme} & N = 30 & N = 40 & N = 50 & N = 60\\ \
    &  $\bar{\epsilon}_{Y_0}$ (Std. Dev.) & $\bar{\epsilon}_{Y_0}$ (Std. Dev.) & $\bar{\epsilon}_{Y_0}$ (Std. Dev.) & $\bar{\epsilon}_{Y_0}$ (Std. Dev.)\\ \hline
   DBSDE & 2.15e-1 (4.19e-3) & 1.83e-1 (5.07e-2) & 1.59e-1 (2.65e-3) & 1.49e-1 (4.12e-3) \\ \hline
   LDBSDE & 3.99e-1 (2.18e-2) & 4.04e-1 (1.60e-2) & 4.21e-1 (1.93e-2) & 4.20e-1 (1.01e-2) \\ \hline
   LaDBSDE & 1.59e-1 (2.78e-2) & 1.69e-1 (2.30e-2) & 1.96e-1 (2.44e-2) & 1.95e-1 (1.43e-2) \\ \hline
   \end{tabular}
  \end{center}
 \caption{The mean absolute errors of $Y_0$ for Example~\ref{ex4} using $d = 100$. The standard deviation is given in parenthesis.}
\label{tab10} 
}
\end{table}

\section{Conclusion}
\label{sec5}
In this work we have proposed the LaDBSDE scheme as a forward deep learning algorithm to solve high dimensional nonlinear BSDEs. It approximates the solution and its gradient based on a global minimization of a novel loss function, which uses local losses defined at each time step including the terminal condition. Our new formulation is achieved by iterating the Euler discretization of time integrals with the terminal condition. The numerical results shows that the proposed scheme LaDBSDE outperforms the existing forward deep learning schemes~\cite{weinan2017deep,raissi2018forward} in the sense of that it does not get stuck in a poor local minima and provide a good approximation of the solution for the whole time domain.

\bibliography{bibfile}
\bibliographystyle{apalike}
\end{document}